\documentclass[fleqn,10pt]{wlscirep}
\usepackage[utf8]{inputenc}
\usepackage[T1]{fontenc}
\usepackage{amsmath}
\usepackage{here}
\usepackage{color}

\newcommand\rev[1]{\textcolor{black}{#1}}

\DeclareMathOperator*{\argmin}{arg\,min}
\newcommand{\bbN}{{\mathbb N}}
\newcommand{\bbR}{{\mathbb R}}
\newcommand{\calN}{{\mathcal N}}

\newcommand\reva[1]{\textcolor{black}{#1}}

\title{Traffic Signal Optimization on a Square Lattice \rev{with Quantum Annealing}}

\author[1,*]{Daisuke Inoue}
\author[1]{Akihisa Okada}
\author[1]{Tadayoshi Matsumori}
\author[2,3]{Kazuyuki Aihara}
\author[1]{Hiroaki Yoshida}
\affil[1]{Toyota Central R\&D Labs., Inc., Bunkyo-ku, Tokyo 112-0004, Japan}
\affil[2]{Institute of Industrial Science, The University of Tokyo, Meguro-ku, Tokyo 153-8505, Japan}
\affil[3]{International Research Center for Neurointelligence, The University of Tokyo, Bunkyo-ku, Tokyo 113-0033, Japan}

\affil[*]{daisuke-inoue@mosk.tytlabs.co.jp}

\begin{abstract}
    The spread of intelligent transportation systems in urban cities has caused heavy computational loads, requiring a novel architecture for managing large-scale traffic.
    In this study, we develop a method for globally controlling traffic signals arranged on a square lattice by means of a quantum annealing machine, namely the D-Wave quantum annealer.
    We first formulate a signal optimization problem that minimizes the imbalance of traffic flows in two orthogonal directions.
    Then we reformulate this problem as an Ising Hamiltonian, which is compatible with quantum annealers.
    The new control method is compared with a conventional local control method for a large 50-by-50 city, and the results exhibit the superiority of our global control method in suppressing traffic imbalance over wide parameter ranges.
    Furthermore, the solutions to the global control method obtained with the quantum annealing machine are better than those obtained with conventional simulated annealing.
    In addition, we prove analytically that the local and the global control methods converge at the limit where cars have equal probabilities for turning and going straight. These results are verified with numerical experiments.
\end{abstract}

\begin{document}

\flushbottom
\maketitle
\thispagestyle{empty}

\section*{Introduction}

For the last two decades, intelligent and efficient transportation systems have been developing, and therefore, control methods for cooperative management of such systems have become increasingly important~\cite{Zhang2011DataDriven,Bishop2005Intelligent,Cheng2015D2D}.
In particular, the adaptive traffic signal operation reflecting the traffic conditions is crucial for avoiding stagnation of traffic flows~\cite{Papageorgiou2003Review,Wei2019Survey}.
\reva{Various methods, which employ techniques such as genetic algorithm~\cite{Gokulan2010Distributed}, swarm intelligence~\cite{Garcia-Nieto2012Swarm}, neural networks~\cite{Srinivasan2006Neural}, and reinforcement learning~\cite{Arel2010Reinforcement,Nishi2018Traffic}, have been proposed for such adaptive control~\cite{Hunt1981Scoot,Roess2004Traffic,Koonce2008Traffic,Faouzi2011Data,Khamis2012Multiobjective,Varaiya2013Maxpressure}}.
In these studies, local control, where the state of each signal is determined from neighboring information, is considered, 
which hardly achieves a global optimum for managing the traffic conditions of the entire city.
Solving a large-scale combinatorial optimization, however, is necessary in order to achieve such a global optimum.
\reva{The difficulty of finding an optimal solution of the latter scales exponentially with the size of the city, because of the computational complexity of the combinatorial optimization.}

\reva{
  Similar computational difficulty frequently appears in other fields. 
  Accordingly, in recent years, various dedicated algorithms and hardware have been developed for solving this issue~\cite{Blum2003Metaheuristics,Puchinger2005Combining,Chakroun2013Combining}. 
  Their main strategy is to focus on solving particular combinatorial optimization problems, which can be transformed into an Ising model. 
  Examples of the specialized hardware include the Coherent Ising Machine provided by NTT Corporation~\cite{Inagaki2016coherent,Hamerly2019Experimental}, 
 the Simulated Bifurcation Machine by Toshiba Corporation~\cite{Goto2019Combinatorial}, 
 and the Digital Annealer by Fujitsu Corporation~\cite{Matsubara2018IsingModel,Aramon2019PhysicsInspired}.
  Among them, the Quantum Annealer 2000Q from D-Wave Systems Inc. has been attracting much attention for its being the world's first hardware implementation of \emph{quantum annealing}~\cite{Kadowaki1998Quantum} using a quantum processor unit.
  In the quantum annealing, a phenomenon called quantum fluctuation is used to simultaneously search candidate solutions of the given problem, which is expected to enable fast and accurate solution search compared with other heuristic search methods~\cite{Kadowaki1998Quantum,Johnson2011Quantum}.
  In this paper, we refer to the method using the 2000Q as the \emph{quantum annealing}.
  Although the quantum annealing is expected to be an effective prescription for the large-scale combinatorial optimization problems, it is not a panacea because the advantage over the classical simulated annealing methods is reduced depending on types of the transformed Ising model.
  Besides the hardware constraints hinder the number of available variables and the class of solvable problems~\cite{Das2008Colloquium}.
  Hence, the search for compatible applications which exploit the quantum annealing power is becoming an active research area~\cite{King2015Benchmarking,McGeoch2013Experimental,Venturelli2016Quantum,OMalley2018Nonnegative,Ohzeki2018Optimization,Inoue2020Model,Ayanzadeh2020Reinforcementa}.
}

In this paper, we propose a method for globally controlling traffic signals in an urban city using the \reva{quantum annealing}.
We consider a situation in which many cars moving on a lattice network are controlled via traffic signals installed at each intersection.
To analytically handle this network,
we consider a simplified situation in which two states are assumed for each signal: traffic is allowed in either the north-south direction or the east-west direction.
The cars moving on the lattice are assumed to choose whether to make a turn or to go straight at an intersection with a given probability.
We then formulate the signal operation problem as a combinatorial optimization problem.
The objective function of the formulated problem is shown to be formally consistent with the Hamiltonian of the Ising model.
The Ising model is a statistical ferromagnetism physics model that represents the behavior of a spin system, and it captures the relation between the microscopic state of spins and the macroscopic phenomena of magnetic phase transitions~\cite{Yang1952Spontaneous,McCoy2014twodimensional,Binder1981Finite,Glauber1963Time}.
\rev{Importantly, the problem reformulated by means of the Ising model, with the aid of a graph embedding technique, is compatible with the class of problems that the 2000Q accepts;}
hence, one can apply quantum annealing to solve the signal optimization problem. 

By reformulating the problem using Ising minimization, this study makes three contributions to signal optimization. 
First, by performing numerical experiments, we confirm the engineering effectiveness of the proposed method using quantum annealing.
Results of experiments using a large city consisting of $50 \times 50$ intersections show that the proposed method achieves \reva{high quality signal operation}, compared with the results of a conventional local control method~\cite{Suzuki2013Chaotic}.
The reformulated optimization problem is also solved using a classical simulated annealing method, but the quantum annealing method is found to give a better solution in a specific parameter domain.
Second, a theoretical correspondence between local and global control methods is found.
We analytically show that the conventional local control is consistent with the solution of the global signal optimization problem at the limit 
where the probability of cars going straight is equal to the probability of them turning.
This result provides a theoretical basis for the numerical prediction of a previous study~\cite{Suzuki2013Chaotic}, where the local control is found to cause phase transitions similar to those of the Ising model. 
The last contribution is the knowledge gained for the cooperative operation of traffic signals.
Our numerical experiments show a strong correlation between a signal and its neighboring signals.
In addition, a strong temporal correlation of signals emerges, that is, the signal display at a certain time is correlated with the displays in the previous several steps.
This spatio-temporal correlation becomes stronger as the straight driving probability of cars increases.
Our results suggest the necessity of signal cooperation for smooth traffic flow, with variation of cooperation strength depending on the rate at which vehicles drive straight.

\section*{Results}\label{sec:formulation}

\subsection*{Traffic Signal Optimization Problem}

Consider $L\times L\ (L\in\bbN)$  roads arranged in east-west and north-south directions with a periodic boundary condition. 
Each road consists of two lanes, one in each direction. 
Traffic signals are located at each intersection to control the flow of vehicles traveling on the roads.
The signal at each node $i$ has one of two states: $\sigma_i=+1$, which allows vehicle flow only in the north-south direction, and $\sigma_i=-1$, which allows vehicle flow only in the east-west direction.
Each car goes straight through each intersection at fixed probability $a\in[0,1]$ and otherwise turns to the left or right with equal probabilities, that is, $(1-a)/2$ for each direction. Figure~1 illustrates this situation.

\begin{figure}[t]
    \centering
    \includegraphics[width=160mm,bb=0.000000 0.000000 825.000000 299.000000]{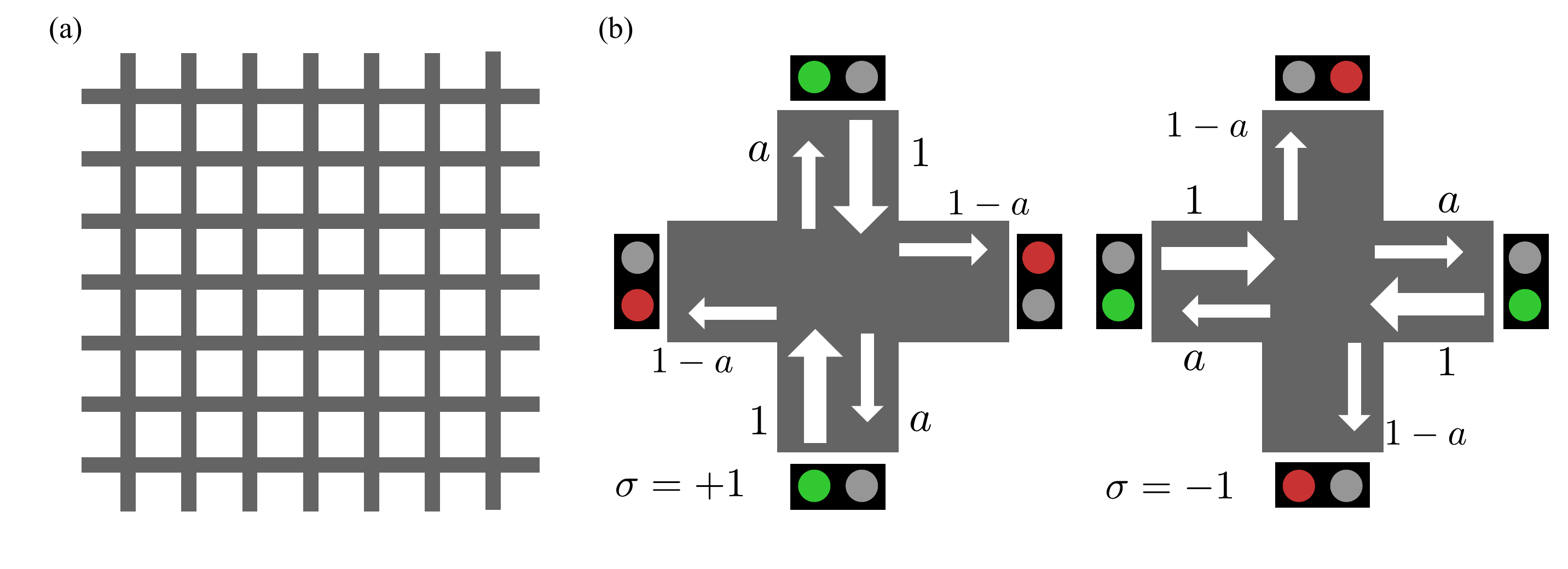}
    \caption{Traffic signal model. (a) Grid pattern of roads. (b) The two states of traffic signals at each intersection.
    In the case of $\sigma=+1$, the vehicles coming from the horizontal direction stop, and the vehicles coming from the vertical direction go straight at the rate of $a$, turn right at the rate of $(1-a)/2$, and turn left at the rate of $(1-a)/2$.
    The rate $1-a$ shown for the horizontal direction is the sum of the vehicles from the two vertical directions.
    In the case of $\sigma=-1$, the roles of the vertical and horizontal directions are reversed. 
    This problem setting is basically following Ref.~\citeonline{Suzuki2013Chaotic}. 
    }
    \label{fig:picture}
\end{figure}

Reference~\citeonline{Suzuki2013Chaotic} shows that the number of vehicles $q_{ij}\in\bbR_+$ in the traffic lane from intersection $j$ to $i$ evolves according to the following difference equation:
\begin{align}\label{eq:traffic_dynamics}
    q_{ij}(t+1) = q_{ij}(t) + \frac{s_{ij}}{2}(-\sigma_i + \alpha \sigma_j),
\end{align}
where $\alpha:=2a-1$, and $s_{ij}\in\{\pm 1\}$ is the direction of the lane from node $j$ to $i$; here, $s_{ij}=+1$ denotes north-south and $s_{ij}=-1$ denotes east-west.
We note here that $q_{ij}$ is normalized by the numbers of cars passing per unit of time.
Precisely, in terms of the mean flux of moving cars $Q_{\mathrm{av}}$ and the dimensional time unit $\Delta t$, 
$t=t^*/\Delta t$ and $q_{ij}=q^*_{ij}/(Q_{\mathrm{av}}\Delta t)$, where $t^*$ is the dimensional time and $q^*_{ij}$ is the number of vehicles in a lane.
\rev{See \emph{Methods} for the detailed derivation of Eq.~\eqref{eq:traffic_dynamics}.}
We define a quantity that represents the deviation of the north-south flow and the east-west flow at each intersection $i$ as
\begin{align}\label{eq:flow}
    x_i(t) := \sum_{j\in \calN(i)} \frac{s_{ij}q_{ij}(t)}{2},
\end{align}
where $\calN (i)$ represents the index of the four intersections adjacent to intersection $i$.
Equation~\eqref{eq:flow} transforms Eq.~\eqref{eq:traffic_dynamics} into a time evolution equation for the flow bias $x(t)$ as follows:
\begin{align}\label{eq:dynamics-merge}
    \mathbf{x}(t+1) = \mathbf{x}(t) + \left( -I + \frac{\alpha}{4} A \right)\pmb{\sigma}(t),
\end{align}
where the flow bias vector is defined as $\mathbf{x}:=[x_1, \ldots , x_{\rev{L^2}}]^\top$ and the signal state vector is defined as $\pmb{\sigma}:=[\sigma_1, \ldots , \sigma_{\rev{L^2}}]^\top$.
The matrix $A\in\bbR^{L^2\times L^2}$ is the adjacent matrix of the periodic lattice graph.

Next, we define the following objective function to evaluate traffic conditions at each time step:
\begin{align}\label{eq:eval_func}
    H(\pmb{\sigma}(t)) := \mathbf{x}(t+1)^\top \mathbf{x}(t+1) + \eta (\pmb{\sigma}(t) - \pmb{\sigma}(t-1))^\top (\pmb{\sigma}(t) - \pmb{\sigma}(t-1)),
\end{align}
where the first term on the right-hand side suppresses the flow bias during the next time step at each intersection, the second term prevents the traffic signal state at each intersection from switching too frequently, and
$\eta\in\bbR_+$ is a weight parameter for determining the ratio of the two terms.
The traffic signal state $\sigma_i(t)$ at each time step is determined so that the objective function \eqref{eq:eval_func} is minimized; that is, we want to find the value of $\pmb{\bar\sigma}(t)$ that satisfies 
\begin{align}\label{eq:optimal_control_problem}
    \pmb{\bar\sigma}(t) =
     \argmin_{\pmb{\sigma}\in\{\pm 1\}^{\rev{L^2}}} H(\pmb{\sigma}(t)).
\end{align}

\subsection*{Ising Formulation and Optimization}
Substituting Eq.~\eqref{eq:dynamics-merge} into Eq.~\eqref{eq:eval_func} gives the following representation:
\begin{align}\label{eq:objective_function}
    \begin{split}
    H(\pmb{\sigma}(t))
        &= \left( \mathbf{x}(t) + \left( -I + \frac{\alpha}{4}A \right) \pmb{\sigma}(t)\right)^\top\left( \mathbf{x}(t) + \left( -I + \frac{\alpha}{4}A \right) \pmb{\sigma}(t)\right)\\
        &\quad + \eta (\pmb{\sigma}(t) - \pmb{\sigma}(t-1))^\top (\pmb{\sigma}(t) - \pmb{\sigma}(t-1))
    \end{split}\\
    \begin{split}
        &= \pmb{\sigma}(t)^\top \left( \left( -I + \frac{\alpha}{4}A \right)^\top \left( -I + \frac{\alpha}{4}A \right) + \eta I \right)\pmb{\sigma}(t)\\
        &\quad + \left( 2\mathbf{x}(t)^\top\left( -I + \frac{\alpha}{4}A \right) - 2\eta \pmb{\sigma}(t-1)^\top \right)\pmb{\sigma}(t) + c(t),
    \end{split}
\end{align}
where $c(t)$ is a constant term that does \emph{not} include $\pmb{\sigma}(t)$.
By defining the variables
\begin{align}
    J &:= \left( -I + \frac{\alpha}{4}A \right)^\top \left( -I + \frac{\alpha}{4}A \right) + \eta I ,\label{eq:J_ij}\\
    h &:=  2\mathbf{x}(t)^\top\left( -I + \frac{\alpha}{4}A \right) - 2\eta \pmb{\sigma}(t-1)^\top ,
\end{align}
we can represent the objective function \eqref{eq:objective_function} as follows:
\begin{align}\label{eq:ising_model}
    H(\pmb{\sigma}(t))= \pmb{\sigma}(t)^\top J \pmb{\sigma}(t) + h \pmb{\sigma}(t) + c(t).
\end{align}
Equation~\eqref{eq:ising_model} is a quadratic form with variables $\{\pm 1\}$, which matches the Hamiltonian form of the Ising model~\cite{Yang1952Spontaneous}.
Hence, solving the signal optimization problem of the objective function \eqref{eq:eval_func} is regarded as equivalent to the problem of finding the spin direction $\sigma_i\in\{\pm 1\}$ that minimizes the Ising Hamiltonian of Eq.~\eqref{eq:ising_model}.
Because the Ising Hamiltonian is fully compatible with the class of problems that the 2000Q accepts, quantum annealing can be applied to solve the signal optimization problem.

We use a city consisting of $50 \times 50$ intersections to consider the signal operation problem, and we compare the results of numerical experiments on the following three methods for traffic control:

\begin{itemize}
    \item Local control, which determines the signal display at each time step with the following local rules:
    \begin{align}
    \begin{cases}\label{eq:suzuki_method}
        \sigma_i(t) \leftarrow +1 & \text{if } x_i(t) \ge +\theta,\\
        \sigma_i(t) \leftarrow -1 & \text{if } x_i(t) \le -\theta.
    \end{cases}
    \end{align}
    Equation~\eqref{eq:suzuki_method} switches the display of the traffic signals to reduce the flow bias when 
    the magnitude of the bias becomes larger than the threshold value $\theta\in\bbR_+$ at each intersection.
    To compare the local control with the optimal control, the value of the switching parameter $ \theta $ is determined such that the common objective function \eqref{eq:eval_func} is minimized.
    For details, refer to \emph{Methods}.
    \item Optimal control with simulated annealing, which reduces Eq.~\eqref{eq:ising_model} at each time step using \emph{simulated annealing}.
    Simulated annealing is an algorithm for finding a solution by examining the vicinity of the current solution at each step and probabilistically determining whether it should stay in the current state or switch to a vicinity state.     
    See Ref.~\citeonline{Suman2006survey} for details of simulated annealing.
    We used the \emph{neal} library provided by D-Wave for executing this algorithm.
    \item Optimal control with quantum annealing, which reduces Eq.~\eqref{eq:ising_model} by using quantum annealing with the D-Wave 2000Q.
    Because the problem size exceeds the size of problems that 2000Q can solve, it is subdivided by the \emph{graph partitioning} technique. 
    We used the \emph{ocean} library provided by D-Wave for executing this algorithm.
    See \emph{Methods} for the detailed procedure.
\end{itemize}

\begin{figure}[t]
    \centering
    \includegraphics[width=160mm,bb=0.000000 0.000000 1879.000000 648.000000]{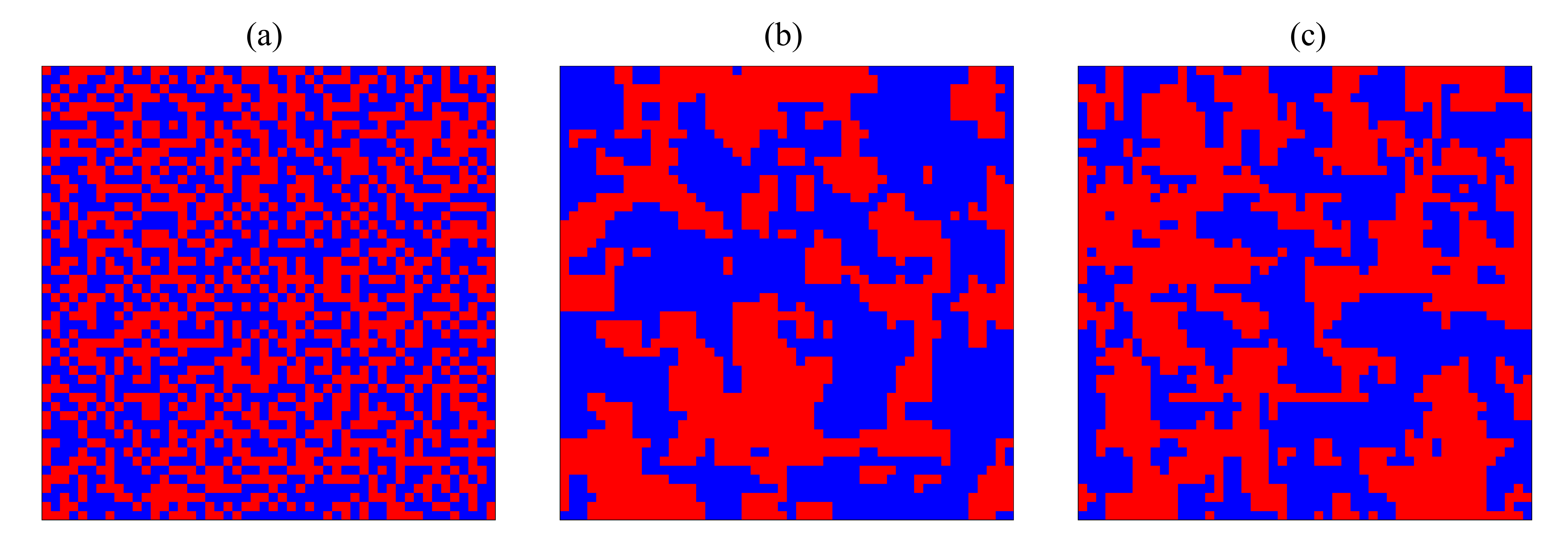}
    \caption{Snapshots of traffic signals under different control methods. (a) Local controller using Eq.~\eqref{eq:suzuki_method}, (b) Global controller optimizing Eq.~\eqref{eq:ising_model} with simulated annealing, and (c) Global controller optimizing Eq.~\eqref{eq:ising_model} with the D-Wave 2000Q.
    Red and blue dots represent vertical and horizontal directions allowed at each crossing, respectively.
    Parameters $\alpha, \eta$, and $L$ are fixed as $\alpha=0.8, \eta=1.0$, and $L=50$, respectively.
    For the D-Wave method, the Hamiltonian is divided into 42 groups and the optimization problem is solved in parallel.
    See \emph{Methods} for details. \reva{(The data are plotted with software Python/matplotlib.)}}
    \label{fig:snapshot50}
\end{figure}

Figure~\ref{fig:snapshot50} shows snapshots of the signal display at time $t = 100$ for $\alpha = 0.8$ and $\eta=1.0$, where $\alpha$ is the parameter related to vehicle's straight driving probability and $\eta$ is the weight parameter in the objective function \eqref{eq:eval_func}.
The flow bias distribution at the initial time $\mathbf{x}(0)$ are generated as random numbers following a uniform distribution of $[-5.0,\ 5.0]$, and the signal states at the initial time $\pmb{\sigma}(0)$ are generated as random numbers following a binomial distribution of $\{\pm 1\}$.
In Fig.~\ref{fig:snapshot50}, blue dots mean that the cars are allowed to pass in the east-west direction, and red dots mean that the cars are allowed to pass in the north-south direction. We observe the synchronization of proximity signals under optimal control [see Figs.~\ref{fig:snapshot50}(b, c)], while the two direction states are distributed rather uniformly under local control [see Fig.~\ref{fig:snapshot50}(a)]. 
\rev{The correlation of proximity signal states is quantitatively analyzed in \emph{Discussion}}.

\begin{figure}[t]
    \centering
    \includegraphics[width=160mm,bb=0.000000 0.000000 825.000000 299.000000]{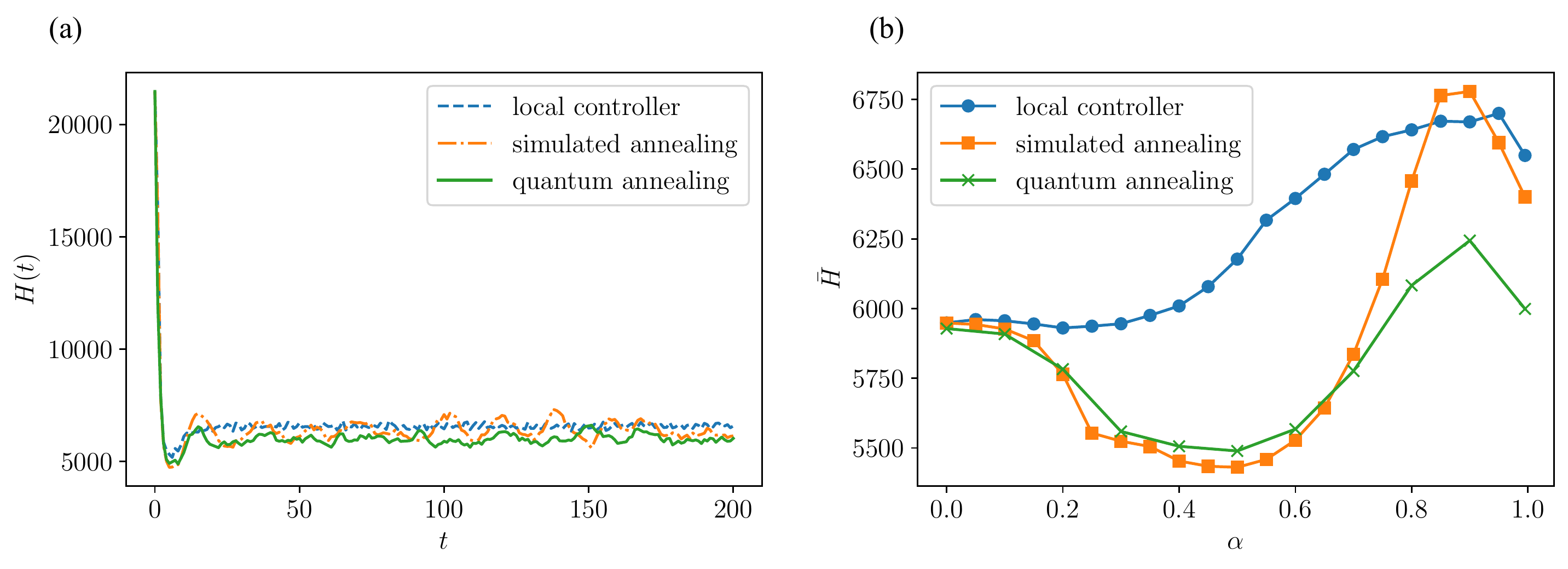}
    \caption{Hamiltonian of Eq.~\eqref{eq:ising_model} under different control methods.
    (a) Time evolution of the Hamiltonian, where the parameters $\alpha, \eta$, and $L$ are fixed as $\alpha=0.8, \eta=1.0$, and $L=50$, respectively.
    (b) Time average of Hamiltonian as functions of $\alpha$, where the parameters $\eta$ and $L$ are the same as those in (a).
    }
    \label{fig:hamiltonian}
\end{figure}

Figure~\ref{fig:hamiltonian}(a) plots the time evolution of the Hamiltonian of Eq.~\eqref{eq:ising_model} for each method when $\alpha = 0.8$ and $\eta=1.0$.
In all three methods, the signals change rapidly over time to reduce the Hamiltonian.
The value of the Hamiltonian in the steady state is the smallest in the quantum annealing method,
followed by the simulated annealing method, and it is the largest under local control.
That is, optimal control using quantum annealing exhibits the best performance among the these methods.
\rev{An attempt to compare with the exact solution has also been made for the same simulation using \emph{Gurobi} package.
Although the full exact solution for the entire time interval was not feasible in a reasonable time because of the large number of variables (2,500 variables), the Hamiltonian averaged over first three steps showed the same order of accuracy.}
\rev{
In Fig.~\ref{fig:hamiltonian}, the response of quantum annealing and simulated annealing is more oscillatory than that of local control.
This is because the objective function of Eq.~\eqref{eq:eval_func} only contains states up to one step ahead.
An optimal value at one time is not necessarily consistent with the optimal values for long time behavior, resulting in more oscillatory response.
If we use an objective function including more than two steps ahead, the oscillatory phenomenon should be more suppressed,
although the latter makes the formulation more complex, hindering the direct use of quantum annealings.}

We examine the effect of changing the parameter $ \alpha $, the vehicle's straight driving probability, on the Hamiltonian of Eq.~\eqref{eq:ising_model}.
The time average of the Hamiltonian of Eq.~\eqref{eq:ising_model}, denoted as $\bar H$, is plotted in Fig.~\ref{fig:hamiltonian}(b).
As $\alpha$ approaches zero, the values of the Hamiltonian for the local and optimal control methods converge to a common value.
This suggests that local control gives the solution to the signal optimization problem at the limit of $\alpha \to 0$. The validity of this conjecture is explored in \emph{Discussion}.
In the interval of $\alpha\in[0.2,\ 0.8]$, the Hamiltonian under optimal control is smaller than that under local control, showing  that the optimum control method exhibits performance better than that of local control in this range.
However, in the simulated annealing method at $ \alpha> 0.8 $, the value of the Hamiltonian is larger than that under the local control method, suggesting that simulated annealing does \emph{not} reach the global optimal solution.
Conversely, under the quantum annealing method, the value of the Hamiltonian is smaller than that under the other two methods, which means that the solution is closer to the global optimum.
\rev{Here, we briefly discuss the slightly better value of $\bar H$ for simulated annealing in a parameter domain of $\alpha\in[0.2,\ 0.8]$, than that for quantum annealing.
In the range of large values of $\alpha$, obtaining an exact solution of Eq.~(4) is hard because of the high impact of the quadratic term.
Actually in this parameter range, the quantum annealing gives better optimization results than the simulated annealing.
On the other hand, regardless of problem to be solved, the quantum annealing generally contains stochastic fluctuations in the solutions~\cite{Ohzeki2018Optimization, King2015Benchmarking}.
When the parameter $\alpha$ is in the intermediate range where the difficulty inherent in the optimization problem is moderate,
both the simulated annealing and the quantum annealing give \reva{high quality solutions}, but the simulated annealing gives slightly better solutions than the quantum annealing because the relative strength of stochastic fluctuations is large.}

\section*{Discussion}\label{sec:analysis}

\begin{figure}[t]
    \centering
    \includegraphics[width=160mm,bb=0.000000 0.000000 825.000000 299.000000]{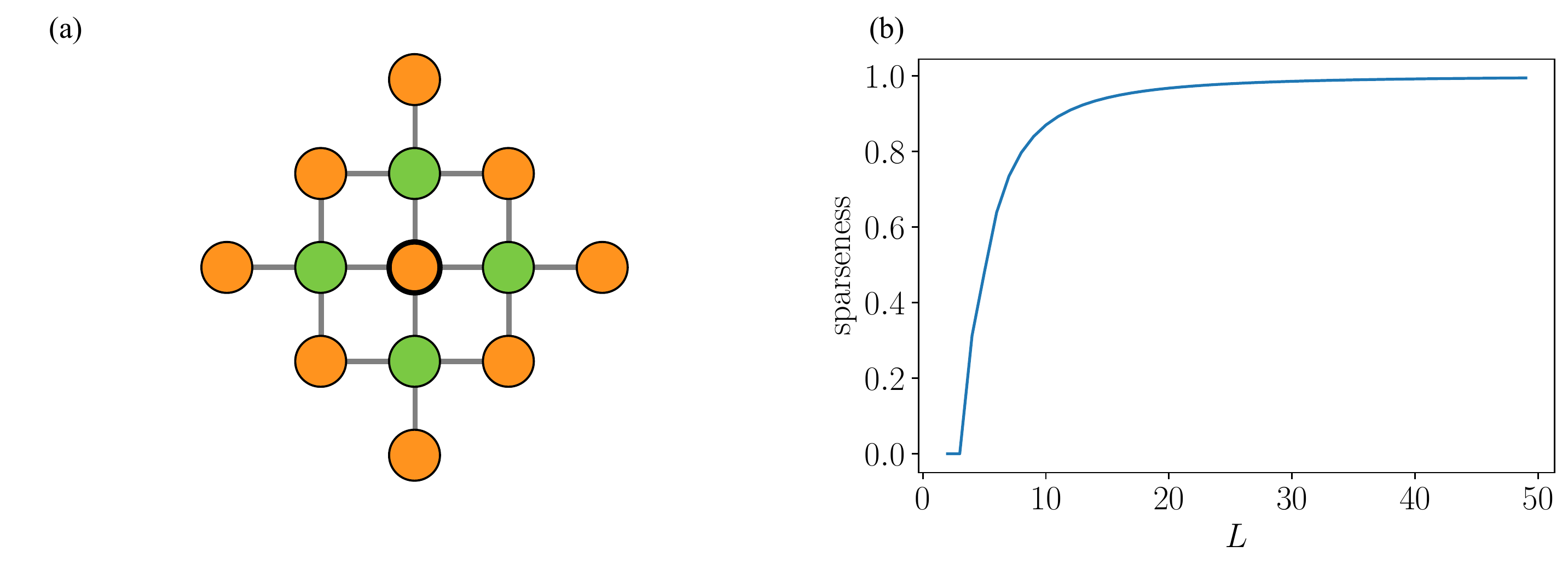}
    \caption{Sparseness of the matrix $J$ in Eq.~\eqref{eq:J_ij}.
    (a) Nodes neighboring the reference node (green) and two nodes away from the reference node (orange) in a lattice graph.
    (b) Sparseness $S_J(L)$ of Eq.~\eqref{eq:sparseness} for different numbers of intersections $L$. 
    }
    \label{fig:sparseness}
\end{figure}

\underline{\it Performance Analysis of Quantum Annealing.} The performance of the D-Wave 2000Q is known to vary depending on the structure of the problem. 
In particular, when the matrix $ J $ in Eq.~\eqref{eq:J_ij} has a sparse structure, the accuracy of the solution is improved~\cite{Hamerly2019Experimental}. 
To check the sparseness of our formulated problem, we examine the value of all components of $J$ in Eq.~\eqref{eq:J_ij}.
First, expanding $ J $ yields the following expression:
\begin{align}\label{eq:J_expand}
    J = (1+\eta)I - \frac{\alpha}{2}A + \frac{\alpha^2}{16}A^\top A, 
\end{align}
where the number of non-zero elements in each column of $ A $ is $ 4 $, because it is equal to the number of degrees of each node in the lattice graph [see the green nodes in Fig.~\ref{fig:sparseness}(a)].
Also, the number of non-zero elements in each column of $ A^\top A $ is $ 9 $ because it coincides with the number of nodes two nodes away from the reference node in the lattice graph [see the orange nodes in Fig.~\ref{fig:sparseness}(a)].
Thus, the number of all non-zero elements in $ J $ is expressed as $ 13L ^ 2 $.
From this, we calculate $S_J(L)$, the sparseness of matrix $J$, defined as the ratio of the number of $0$-valued elements and the number of all elements in the matrix:
\begin{align}\label{eq:sparseness}
    S_J(L) = \frac{L^4-13L^2}{L^4},
\end{align}
where we confirm that $S_J(L)\to 1$ as $L\to\infty$.
In Fig.~\ref{fig:sparseness}(b), we plot $S_J(L)$ given in Eq.~\eqref{eq:sparseness}, to show that the sparseness of matrix $ J $ increases as increasing city size. 
This allows us to expect that the performance of the D-Wave 2000Q is enhanced in the case of the signal optimization problem for a rather large cities, such as $L=50$, the one considered in the present paper.

\begin{figure}[t]
    \centering
    \includegraphics[width=160mm,bb=0.000000 0.000000 825.000000 299.000000]{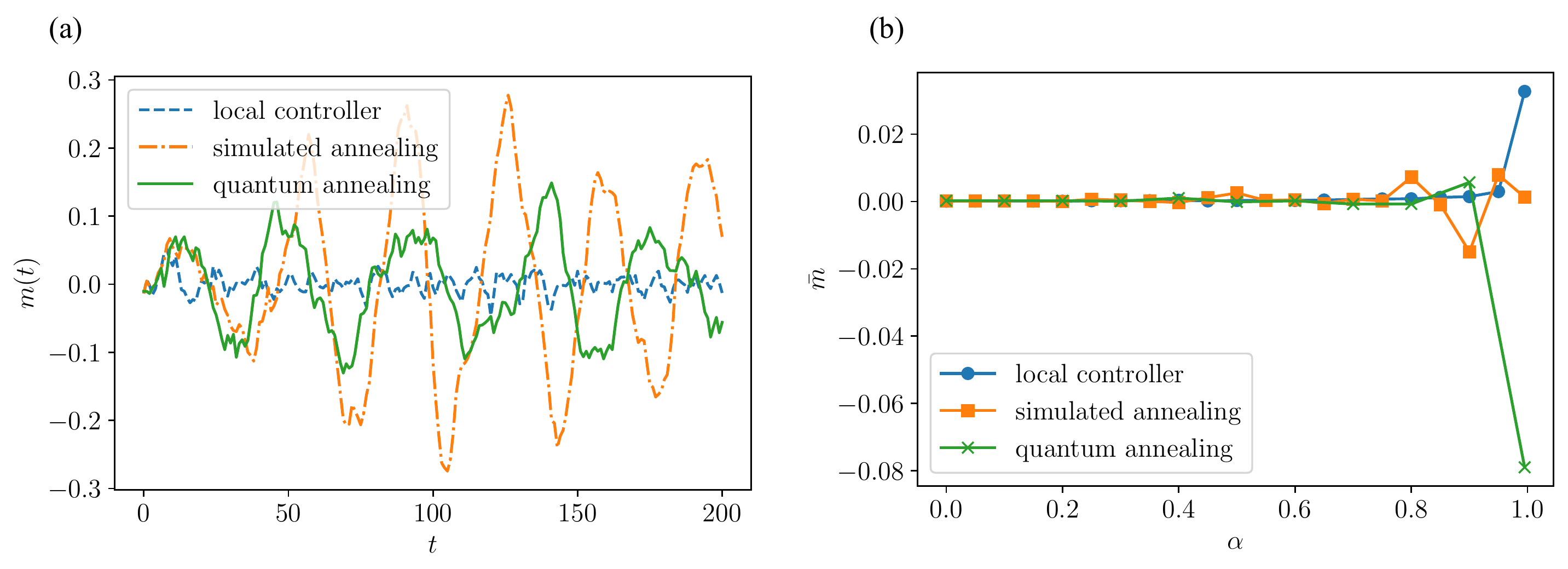}
    \caption{Magnetization of Eq.~\eqref{eq:magnetization} under different control methods.
    (a) Time evolution of magnetization. 
    Parameters $\alpha, \eta$, and $L$ are fixed as $\alpha=0.8, \eta=1.0$, and $L=50$, respectively.
    (b) Time average of magnetization as a function of $\alpha$.  
    Parameters $\eta$ and $L$ are the same as those in (a).
    }
    \label{fig:magnetization}
\end{figure}

\vspace{2mm}
\noindent
\underline{\it Local and Optimal Control Correspondence.} 
As shown in Fig.~\ref{fig:magnetization},
when the parameter $\alpha$ of Eq.~\eqref{eq:traffic_dynamics} is sufficiently small, the local control of Eq.~\eqref{eq:suzuki_method} approaches the optimal control that is the solution of Eq.~\eqref{eq:optimal_control_problem}.
When $\alpha\approx 0$ is valid, the term associated with $ \alpha $ in Eq.~\eqref{eq:ising_model} can be ignored, yielding
\begin{align}
    J &\approx (1+\eta)I,\label{eq:J_approx}\\
    h &\approx  - 2\mathbf{x}(t)^\top - 2\eta \pmb{\sigma}(t-1)^\top.\label{eq:h_approx}
\end{align}
Because $J$ in Eq.~\eqref{eq:J_approx} is a diagonal matrix, the first term $\pmb{\sigma}(t)^\top J \pmb{\sigma}(t)$ on the right-hand side of Eq.~\eqref{eq:ising_model} is a constant that does \emph{not} depend on $\pmb{\sigma}$.
Therefore, the minimizer of $H(\pmb{\sigma}(t))$ is determined depending only on the sign of $h$ in Eq.~\eqref{eq:h_approx}, that is, 
\begin{align}\label{eq:sigma_local_mid}
    \bar\sigma_i(t) =&
    \begin{cases}
        1 & \text{if } x_i(t) + \eta\sigma_i(t-1)\ge 0,\\
        -1 & \text{if } x_i(t) + \eta\sigma_i(t-1) < 0,
    \end{cases}
\end{align}
for all $i=1,\ldots,\rev{L^2}$. 
By transforming Eq.~\eqref{eq:sigma_local_mid}, we obtain
\begin{align}\label{eq:opt_input}
  \bar\sigma_i(t) =& 
  \begin{cases}
  1 & \text{if } x_i(t) \ge \eta, \\
  -1 & \text{if } x_i(t) \le -\eta,\\
  \sigma(t-1) & \text{otherwise},
  \end{cases}
\end{align}
for all $i=1,\ldots,\rev{L^2}$.
The control method of Eq.~\eqref{eq:opt_input} is equivalent to the local control \eqref{eq:suzuki_method} in Ref.~\citeonline{Suzuki2013Chaotic}.

Because $ \alpha = 0 \Leftrightarrow a = 0.5 $ holds, this optimality means that an appropriate vehicle turning rate autonomously eases the flow bias in the local control laws.
In addition, the occurrence of this magnetic transition for the signal display, stated in Ref.~\citeonline{Suzuki2013Chaotic}, is consistent with the fact that local control in Eq.~\eqref{eq:suzuki_method} actually minimizes the Ising Hamiltonian in Eq.~\eqref{eq:ising_model}.
However, note that the optimality of local control is valid only when $ \alpha \approx 0 $, but \emph{not} when $ \alpha \to 1 $, where the phase transition occurs.

\vspace{2mm}
\noindent
\underline{\it Signal Synchronization Analysis.}
To analyze the signal correlation observed in Fig.~\ref{fig:snapshot50}, we calculate the magnetization, which is regarded as an important quantity in the Ising model:
\begin{align}\label{eq:magnetization}
    m(t) := \frac{1}{L^2}\sum_{i=1}^{L^2} \sigma_i(t).
\end{align}
In the Ising model, this value represents the spin bias of the entire system, and it is an indicator of ferromagnetic transitions in the system.
Figure~\ref{fig:magnetization}(a) shows the time variation of magnetization $ m(t) $.
The value of magnetization remains small under local control, whereas it becomes significantly larger under both optimal control methods (simulated annealing and quantum annealing).
For each method,  at $\alpha=0.8$, 
the response of the magnetization oscillates or fluctuates around zero.
To confirm this observation, the time average of the magnetization of Eq.~\eqref{eq:magnetization}, denoted as $ \bar m $, is  plotted in Fig.~\ref{fig:magnetization}(b).
Here, the ferromagnetic transition at $ \alpha \to 1 $, that is, the finite value of $\bar m$,
is observed for the magnetization under local control, which was originally reported in Ref.~\citeonline{Suzuki2013Chaotic}.
Also, under optimal control, the time average of the magnetization $\bar m$ takes a large value when $ \alpha \to 1 $, which shows that a ferromagnetic transition similar to that under local control occurs under optimal control.

\begin{figure}[t]
    \centering
    \includegraphics[width=160mm,bb=0.000000 0.000000 825.000000 299.000000]{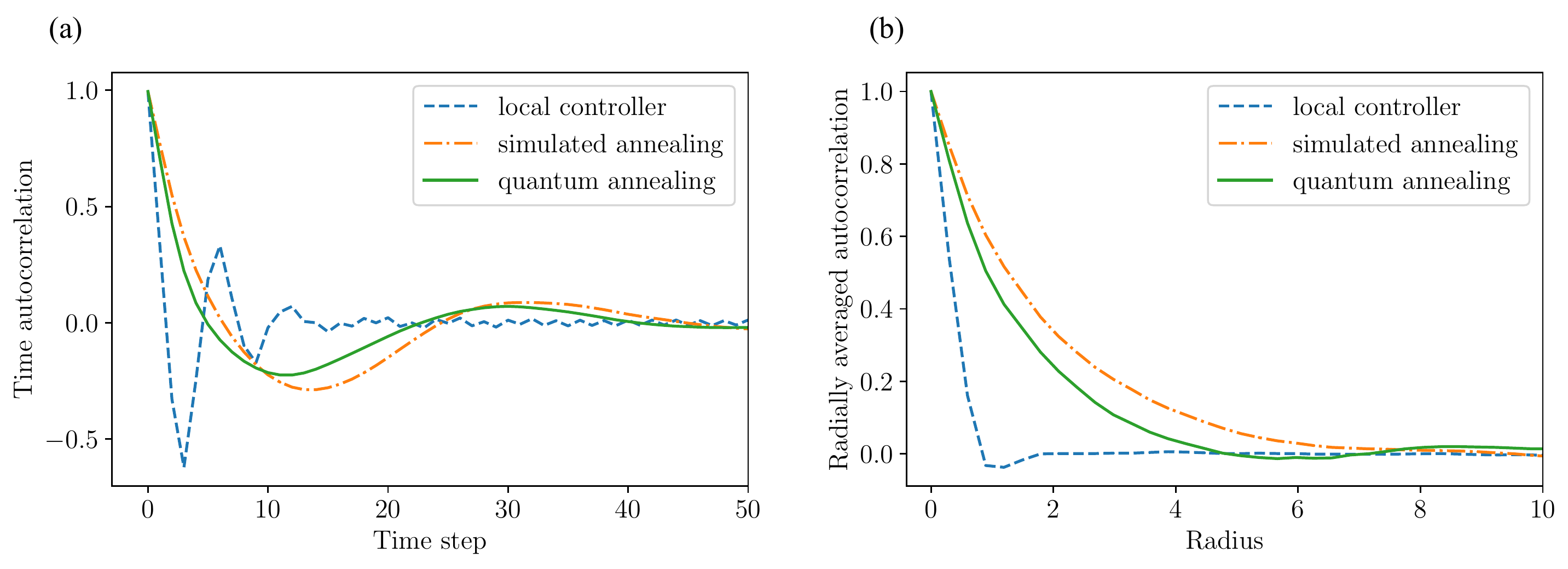}
    \caption{Time and spatial autocorrelation functions for different control methods.
    (a) Time autocorrelation function and (b) Radially averaged spatial autocorrelation function. 
    Parameters $\alpha, \eta$, and $L$ are fixed as $\alpha=0.8, \eta=1.0$, and $L=50$, respectively.}
    \label{fig:autocorrelation}
  \end{figure}

In addition to the ferromagnetic transition, the large amplitudes observed under optimal control
are indeed a quantification of the synchronization of proximity signals observed in Fig.~\ref{fig:snapshot50}.
For further analysis of this synchronization, we also evaluate two types of autocorrelation functions.
Figure~\ref{fig:autocorrelation}(a) shows the autocorrelation function obtained from the time-series data of the signal state $\sigma_i(t)$ for $t\in[0,200]$.
Here, the autocorrelation function is computed at all intersections, and the average value is displayed in Fig.~\ref{fig:autocorrelation}(a).
Under local control, there is a negative correlation peak around $ t = 3 $, which means that the signals switch approximately every $ 3 $ time steps.
In contrast, under optimal control, the negative correlation peak is in the interval of $ t = [10,15] $ steps, and the same state is maintained for a time longer than that under local control.
\rev{In general, excessive signal switching is undesirable from a traffic engineering standpoint, and the optimization-based method overcomes this issue.}
Next, Fig.~\ref{fig:autocorrelation}(b) shows the correlation between the display of signals at one intersection and another intersection, with the distance between the intersections as a parameter.
Here, the correlation function is calculated for all the intersections for fixed time $t=100$, and the average value thereof is plotted.
In Fig.~\ref{fig:autocorrelation}(b), the distance is normalized to make the distance of adjacent intersections equal to $ 1 $.
There is almost no correlation between adjacent signals under local control, while there is a positive correlation of up to 4--6 adjacent intersections under optimal control.

\begin{figure}[t]
  \centering
  \includegraphics[width=160mm,bb=0.000000 0.000000 825.000000 299.000000]{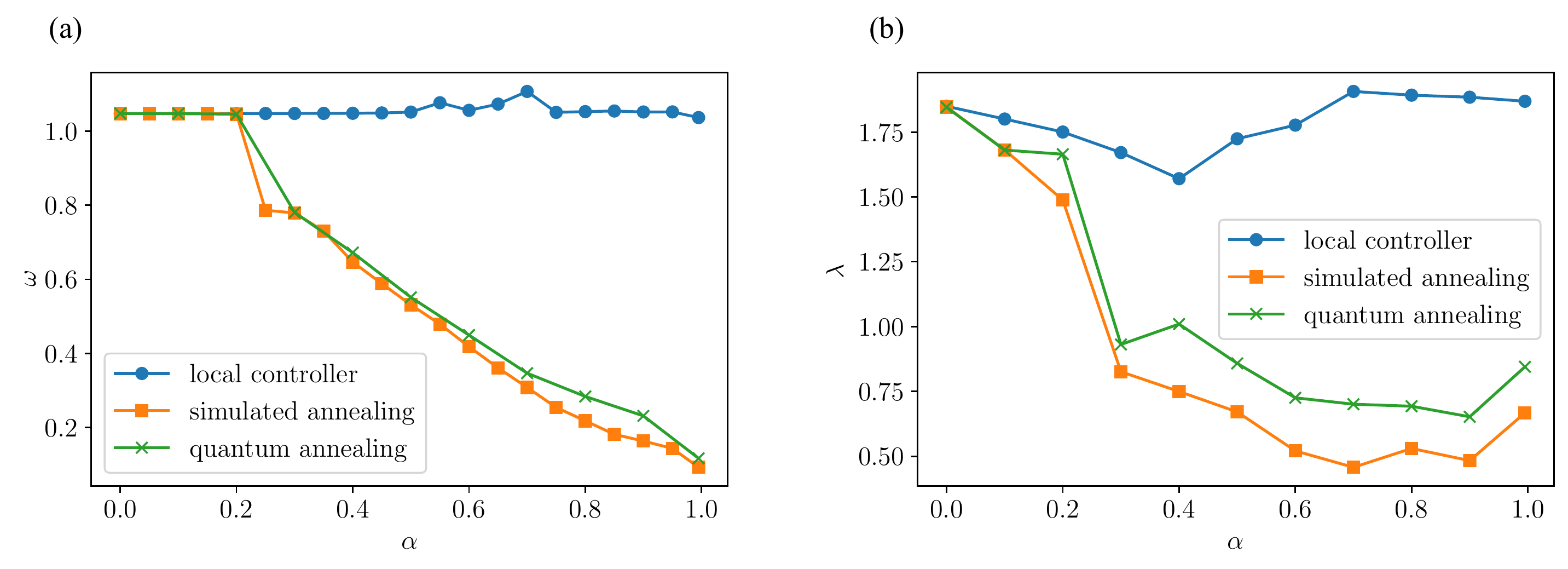}
  \caption{Parameters extracted from time and spatial autocorrelations, as functions of $\alpha$ for different control methods.
  (a) Time autocorrelation function frequency $\omega$ versus $\alpha$ and (b) Radially averaged autocorrelation decay rate $\lambda$ versus $\alpha$. 
  Parameters $\eta$ and $L$ are fixed as $\eta=1.0$ and $L=50$, respectively.}
  \label{fig:autocor_alpha_sweep50}
\end{figure}

Then, we extract quantities from these correlation functions to investigate the effect of $\alpha$.
First, considering that both the temporal and spatial autocorrelations in Fig.~\ref{fig:autocorrelation} decay while oscillating, both functions are fitted with the following equation:
\begin{align}\label{eq:autocor_fitting}
  R(z) = \exp(-\lambda z)\cos(\omega z),
\end{align}
where $\lambda$ represents the damping rate coefficient, $ \omega $ represents the vibration frequency coefficient, and $z\in\bbR_+$ represents different variables, i.e., the time $t$ for the time autocorrelation function and the distance between intersections for the spatial autocorrelation function.
Figure~\ref{fig:autocor_alpha_sweep50}(a) plots $ \omega $ values obtained by fitting Eq.~\eqref{eq:autocor_fitting} to the time autocorrelation, 
as a function of $ \alpha $.
Under local control, the vibration frequency is $ \omega \approx 1 $ regardless of the value of $ \alpha $, 
while $ \omega $ decreases as increasing $ \alpha $ under optimal control.
%
\rev{This suggests that the frequency of signal switching reduces as the vehicle straight driving rate increases in order to guarantee optimality.
In view of the large difference in $\omega$ between the local control and the optimization-based controls
for large values of $\alpha$, 
we expect that optimization-based signal controls are particularly effective
in preventing excessive switching for high vehicle straight driving rates.}
Next, we show in Fig.~\ref{fig:autocor_alpha_sweep50}(b) the value of $ \lambda $ obtained by 
fitting Eq.~\eqref{eq:autocor_fitting} to the spatial autocorrelation, as a function of $\alpha$.
Under local control, the correlation decreases with an attenuation factor of $ \lambda \approx 1.75 $, regardless of the value of $ \alpha $.
In contrast, under optimal control, $ \lambda $ decreases as $ \alpha $ increases, which means that the signal displays between the more distant intersections remain correlated.
These observations show that the synchronization of proximity signals in time and space becomes important for achieving a balanced traffic flow as the probability of vehicles going straight increases.

\vspace{2mm}
\noindent
\rev{\underline{\it Effect of Parameter $\eta$.} }
\begin{figure}[t]
    \centering
    \includegraphics[width=80mm,bb=0.000000 0.000000 393.083371 267.635810]{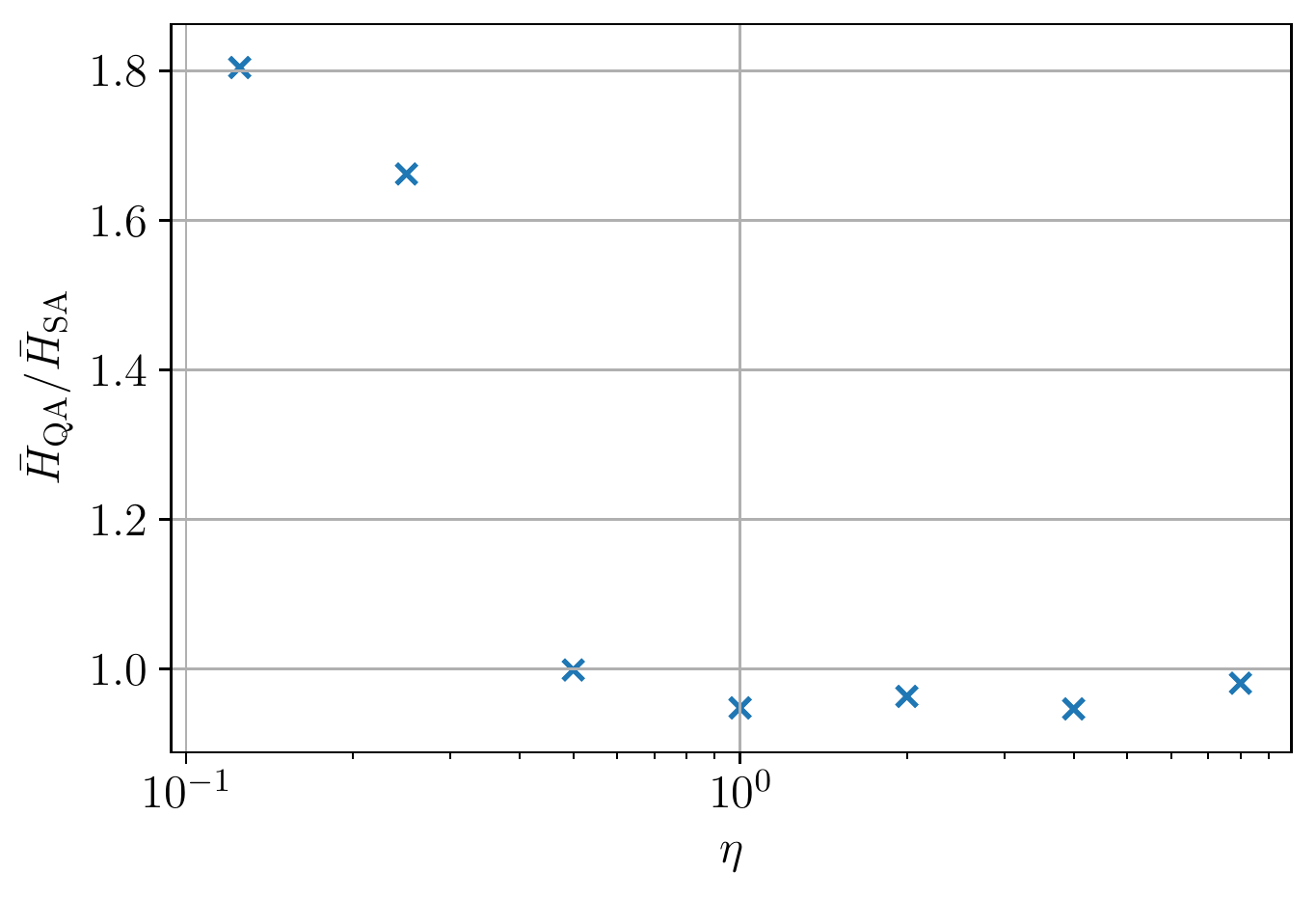}
    \caption{\rev{
    The time averaged values of the Hamiltonian for the simulated annealing ($\bar H_\text{SA}$) and the quantum annealing ($\bar H_\text{QA}$) as a function of $\eta$.
    Parameters $\alpha$ and $L$ are fixed at $\alpha=0.8$ and $L=50$.}}
    \label{fig:eta_change}
  \end{figure}
\rev{
Here we examine the effect of parameter $\eta$, which controls the priority of the smoothness of the entire traffic flow to the signal switching frequency in the objective function of Eq.~\eqref{eq:eval_func}.
The objective function is designed such that the priority is to smooth the flow of the car for small value of $\eta$, and inversely for a large value of $\eta$, preventing excessive signal switching is prior.
The time average of the objective function $\bar H$ is obtained for various values of $\eta$; $\eta\in\{ 0.125,0.25,0.5,1,2,4,8\}$.
We show the results in Fig.~\ref{fig:eta_change} for the simulated annealing ($\bar H_\text{SA}$) and the quantum annealing ($\bar H_\text{QA}$), where the ratios $\bar H_\text{QA}/\bar H_\text{SA}$ is plotted; $\bar H_\text{QA}/\bar H_\text{SA}<1$ means that the quantum annealing is better than the simulated annealing, and vice versa.
The quantum annealing method shows better performances for $\eta$ larger than $0.5$,
and the simulated annealing is better for $\eta$ smaller than $0.5$.
This suggests that the quantum annealing works better when the priority is on preventing excessive signal switching.
}

\vspace{2mm}
\noindent
\reva{\underline{\it Future improvements.} }
\reva{
  Here we discuss three possible improvements of the results obtained in this study.
  First, we expect that the solution is improved by using the most recently released D-Wave's machine.
  The D-Wave machine used in this study has $2,048$ qubits that are connected with \emph{chimera structure}, in which closely connected 8-qubit units are arranged~\cite{Boothby2016Fast}.
  Since the chimera structure is sparser than the fully-connected structure, the representation of arbitrary Ising problems requires a process called \emph{minor-embedding} to map logical variables to physical qubits.
  This process however significantly reduces the number of available qubits, and also deteriorates the computational accuracy.
  Very recently, D-Wave has updated the hardware with a new graph structure called the \emph{pegasus structure}.
  The number of qubits has increased  from 2,048 to 5,024, and the maximum number of connections in the graph structure has increased from 8 to 15~\cite{W.Johnson2018Future}. 
  These improvements allow us to deal with much larger problems, and to realize efficient embedding with sacrificing less physical qubits, than the previous D-Wave machine. 
 For the proposed method, this enhancement will significantly reduce the number of divisions of Hamiltonian (see \emph{Methods} for details), which contributes to fast and high-accurate computations.
}

\reva{
  Next, adjusting the hyper-parameters of the solver could improve the performance of our method.
  The  D-Wave machine contains a few hyper-parameters, such as, \emph{number of samplings}, \emph{chain strength}, and \emph{post processing}.
  We left most of the parameters at their default values because we focus on examining ability of the quantum annealing to solve the traffic signal optimization problem. 
  However, as these hyperparameters affect the optimized result, more careful tuning of these parameters may achieve faster and more accurate calculations.
  An error mitigation scheme proposed in Ref.~\citeonline{Ayanzadeh2020PostQuantum} could enhance the performance.
  We remark here that the problem formulated in this study has the form of Ising model that is solvable by using several dedicated computers other than D-Wave machines~\cite{Inagaki2016coherent,Hamerly2019Experimental,Matsubara2018IsingModel,Goto2019Combinatorial}.
  Since the development of these dedicated machines is expected to further accelerate, the proposed framework for traffic flow control will be more generally available in the future.
}

\reva{
  We finally discuss a further improvement toward application to a real city.
  The parameter $\alpha$ in our model is expected to be relatively large in a city with many rational players, because in a grid-like city, each vehicle can reach its destination from any starting point with only one right or left turn.
  In our experiments, the size of the grid $L$ is empirically determined as $L=50$ so that it would be comparable to the size of typical grid cities in the world (Kyoto, Japan; Barcelona, Spain; La Plata, Argentina, etc.).
  It is however desirable to identify these parameters in advance using real-world data.
  Since the constant probability of each vehicle driving straight ahead and the grid topology of the city are both idealistic assumptions, our traffic signal control method has to be further improved by relaxing these assumptions.
}

\section*{Methods}

\subsection*{\rev{Derivation of Traffic Model}}
\rev{
Here we derive the model shown in Eq.~\eqref{eq:traffic_dynamics}.
Let $q_{ij}^*(t)$ be the numbers of cars that exist between the intersections $i$ and $j$ and $\Delta t$ be the minimum time interval at which a signal is switched.
We denote by $Q_\text{av}$ the average flow rate of cars passing during $\Delta t$.
Then the change in the numbers of cars from time $t^*$ to the next time $t^*+\Delta t$ is represented as
\begin{align}
  q_{ij}^*(t^*+\Delta t) = q_{ij}^*(t^*) + \frac{s_{ij}}{2}(-\sigma_i(t^*)+\alpha\sigma_j(t^*)) Q_\text{av} \Delta t, 
  \quad q_{ij}^*(0) = q_{ij}^{*0}.
 \label{eq:traffic_dimensional}
\end{align}
By normalizing equation \eqref{eq:traffic_dimensional}
with $t:=t^*/\Delta t$ and $q:=q^*/(Q_\text{av} \Delta t)$, we obtain the following equation:
\begin{align}\label{eq:traffic_dynamics_2}
  q_{ij}(t+1) = q_{ij}(t) + \frac{s_{ij}}{2}(-\sigma_i(t) + \alpha \sigma_j(t)),
  \quad q_{ij}(0) = q_{ij}^{*0}/(Q_\text{av} \Delta t),
\end{align}
which is essentially identical to Eq.~\eqref{eq:traffic_dynamics}.
In this paper, we consider the result of solving Eq.~\eqref{eq:traffic_dynamics_2}. 
The dimensional time and the actual numbers of cars are 
recovered with inverse transformation of the above normalization.
}

\rev{
We remark on the numbers of cars, speed, and minimum signal switching interval on the model.
First, a solution of the model in Eq.~\eqref{eq:traffic_dynamics} 
is valid for an arbitrary numbers of cars.
For example, the solution for $\bar q^*(0)=\gamma q^*(0)$ with some $\gamma\in\bbR_+$, 
is obtained by setting $\bar Q_\text{av} = \gamma Q_\text{av}$ because the average flow rate is defined by ``vehicle density'' $\times$ ``average speed''.
Second, let us consider the case of the average speed multiplied by $\gamma$.
The average flow rate $Q_\text{av}$ should then be $\tilde Q_\text{av} = \gamma Q_\text{av}$, while $q^*(0)$ remains the same.
Therefore, while Eq.~\eqref{eq:traffic_dynamics_2} normalized by $\tilde Q_\text{av}$ does not apparently change,
the initial value should be appropriately adjusted $\tilde q_{ij}(0) = q_{ij}^{*0}/(\tilde Q_\text{av} \Delta t) = q_{ij}(0)/\gamma$.
Similarly, for the case of $\Delta \hat t := \gamma \Delta t$,
the normalized equation \eqref{eq:traffic_dynamics_2} does not apparently change, but the initial value should be $\hat q_{ij}(0) = q_{ij}^{*0}/(Q_\text{av} \Delta \hat t) = q_{ij}(0)/\gamma$.
}

\subsection*{Parameter Identification for Objective Function}\label{sec:eta-theta}

As stated in \emph{Discussion}, 
a direct correspondence between the optimal control and local control
is established for small values of $ \alpha $,
with the apparent relation $ \theta = \eta $ 
between the local control switching constant $ \theta $ in Eq.~\eqref{eq:suzuki_method} and the optimal control weight parameter $ \eta $ in Eq.~\eqref{eq:eval_func}.
To make a systematic comparison for an arbitrary value of $\alpha$, 
however, we still need to construct a protocol
to determine the values of $\theta$ and $ \eta $.
The strategy is described as follows.
Given a value of $ \eta $, we select a value of $\theta$, denoted by $\hat \theta$, from a candidate set $ \Theta $
via the following auxiliary numerical analysis:
\begin{enumerate}
    \item For one value of $\theta$ in the set $\Theta$,
    numerical simulation using local control \eqref{eq:suzuki_method} is performed 
    to obtain time series data $\mathbf{x}(t)$ and $\pmb{\sigma}(t)$.
    The value of the objective function \eqref{eq:eval_func} using the given $ \eta $ is calculated from the obtained time series data.
    This time average is denoted as $  \bar H(\theta) $.
    \item Step~1 is performed for all $\theta$ in $\Theta $ to find $ \hat \theta $ that minimizes the time average $ \bar H $, that is, $\hat \theta = \argmin_{\theta\in\Theta}\bar H(\theta)$.
\end{enumerate}

\begin{figure}[t]
    \centering
    \includegraphics[width=160mm,bb=0.000000 0.000000 825.000000 299.000000]{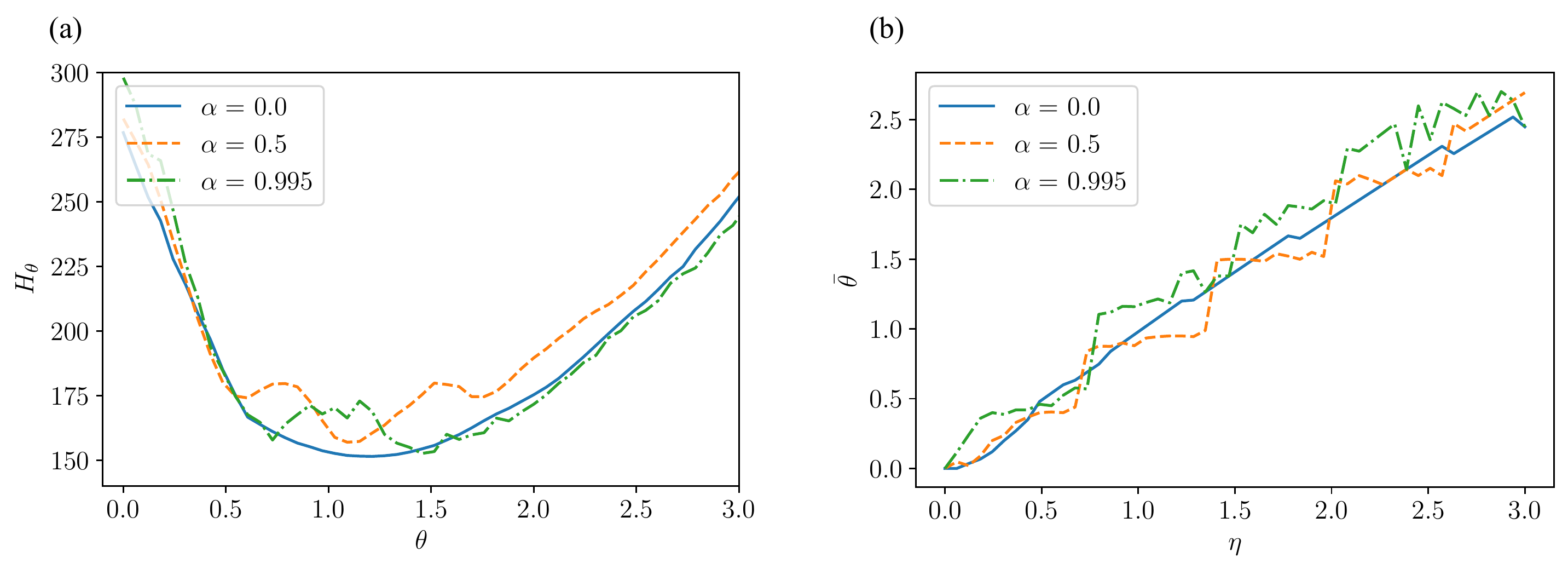}
    \caption{Correspondence between $\eta$ and $\theta$. (a) Time average of the objective function $\bar H$ versus $\theta$, when the value of $\eta$ is fixed as $\eta=1.0$. The cases with $\alpha\in\{0.0,0.5,0.995\}$ are shown. (b) $\hat\theta$ versus $\eta$ for $\alpha\in\{0.0,0.5,0.995\}$.}
    \label{fig:eta-theta}
\end{figure}

We plot the result of the above procedure in Fig.~\ref{fig:eta-theta}.
Figure~\ref{fig:eta-theta}(a) shows $ \bar H $ against $ \theta $ when $ \eta$ is fixed as $ \eta = 1.0 $.
When $ \alpha = 0 $, $ \bar H $ is a convex function and indeed $ \hat \theta \approx \eta $ is satisfied.
However, for larger values of $\alpha$,  $ \bar H $ becomes non-convex, and 
particularly for $\alpha=0.995$, the relation $ \hat \theta = \eta $ no longer holds.
Figure~\ref{fig:eta-theta}(b) shows the value of $ \hat \theta $ that minimizes $H$ versus $ \eta $ for the interval $\eta \in [0.0, \ 3.0] $.
When $ \alpha = 0 $, the linear relation $ \hat \theta = \eta $ approximately holds, but when $ \alpha \ne 0 $, this relation breaks down and some discontinuities appear.
These discontinuities correspond to the changes in the local minima observed in Fig.~\ref{fig:eta-theta}(a).

\subsection*{\rev{Implementation on D-Wave 2000Q}}\label{sec:dwave_soltion}

\reva{
  All experiments in this study are conducted on a Linux computer with 64 GB of memory and a clock speed of 3.70 GHz. 
  All methods are implemented using the programming language Python (version 3.7).
  %
  We use \emph{DW\_2000Q\_VFYC\_5} as a machine solver with the aid of D-Wave's \emph{ocean} library for the actual implementation.
  Here, the VFYC solver partially emulates some qubits that are temporarily unavailable because of hardware failures~\cite{manual}.
  The number of samplings is specified through a parameter named \emph{num\_reads}, which we set $100$ in all experiments.
  The validity of this parameter setting is confirmed by preliminary experiments using several candidate parameters.
  For embedding the logical variables to the physical qubits on the D-Wave machine, we use a tool called \emph{minorminer} (Apache license 2.0), which is a heuristic embedding method in ocean library~\cite{Cai2014practicala}.
  We perform embedding operation every time when the problem is sent to 2000Q in order to average out the bias of the embedding quality.
  For the simulated annealing, \emph{neal} solver in \emph{ocean} library is used.
  This solver also allows us to specify the number of samplings through a parameter called num\_reads, which we set $100$, i.e., the same value as the one in the quantum annealing. 
}

\begin{figure}[t]
    \centering
    \includegraphics[width=80mm,bb=0.000000 0.000000 734.400000 734.400000]{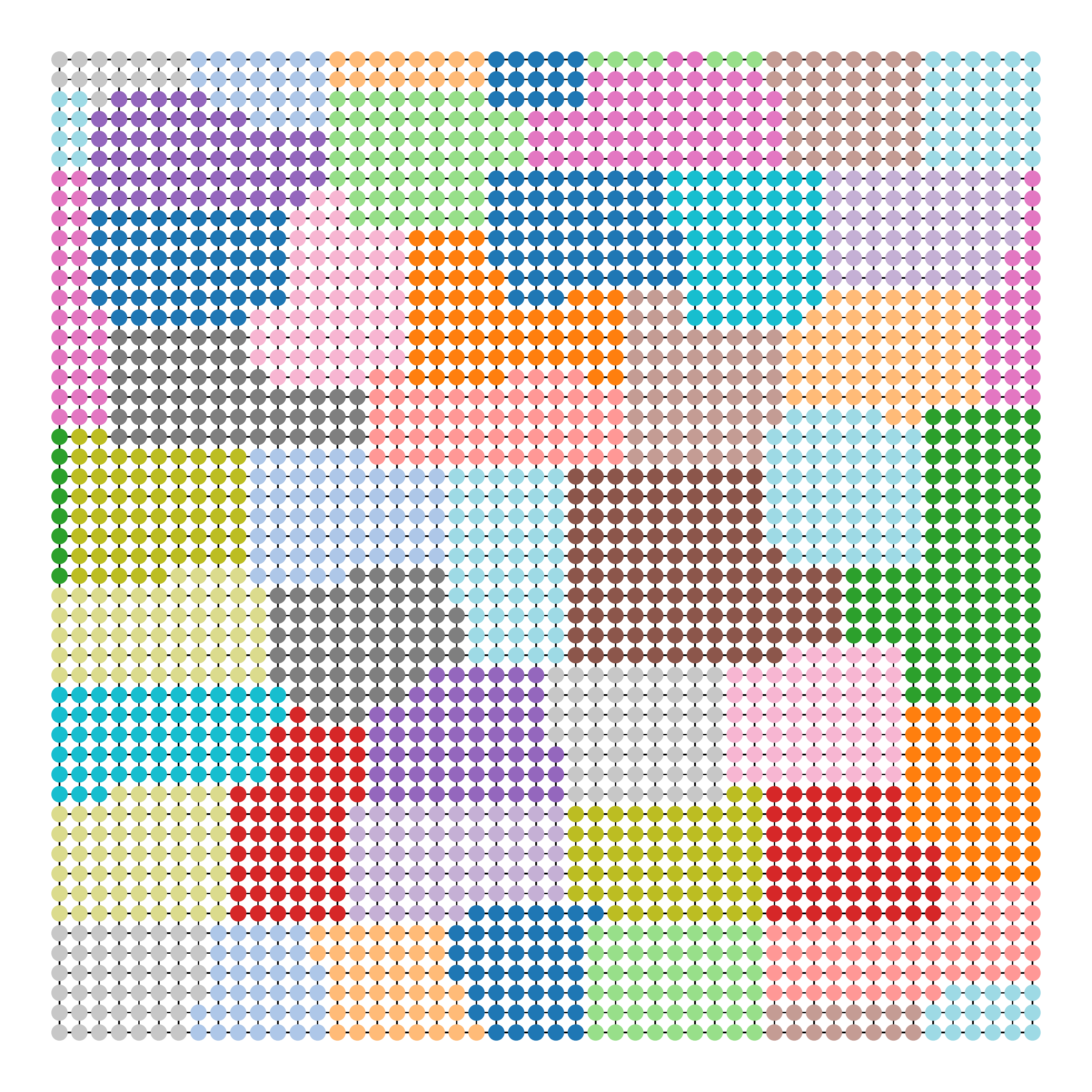}
    \caption{Graph partitioning using Metis. Each node represents a component of the Hamiltonian coefficient matrix $ J $ in Eq.~\eqref{eq:ising_model}, and the color of each node indicates the group to which the component belongs. \reva{(The data obtained using Metis 5.1.0 are plotted with software Python/matplotlib.)}}
    \label{fig:metis}
\end{figure}

For the chimera structure in the 2000Q, $N^2/4$ physical qubits are necessary for embedding $N$-variable problem for the worst case, which means that the maximum number of variables that the 2000Q is capable of handling is as small as 64. 
This implies that $L^2\le 64\Leftrightarrow L\le 8$ must be satisfied for the number of roads $L$ in our problem setting. 
A method exists for solving a problem that exceeds the size limitation: to divide the Hamiltonian variable of Eq.~\eqref{eq:ising_model} into several groups and minimize the approximate Hamiltonian for each group.
We define the traffic signal state vector of the $j$th group as $\pmb{\sigma}^{j} := [\sigma_{i_1}, \sigma_{i_2}, \ldots , \sigma_{i_m}]^\top$, where $ i_1, i_2, \ldots, i_m $ are subscripts of variables included in the $j$th group.
Then, we define the Hamiltonian of the group $j$ as
\begin{align}\label{eq:split_hamiltonian}
    H^j(\pmb{\sigma}^j(t)) := \pmb{\sigma}^{j}(t)^{\top} J_{jj}\pmb{\sigma}^j(t) + 
    (h_j + \pmb{\sigma}^{\bar j}(t)^{\top} J_{\bar j j} )\pmb{\sigma}^j(t),
\end{align}
where $J_{jj}$ is a matrix extracting the $ (j, j) $th components of matrix $ J $ in Eq.~\eqref{eq:ising_model}.
Similarly, $ h_j $ is a vector obtained by extracting the $j$th component of $ h $.
The index $ \bar j $ represents the set of variables \emph{not} belonging to group $j$.
One naive approximation is to regard the variables outside group $ j $ as constant.
This allows the annealing machine to deal with a Hamiltonian exceeding the limitation, but at the same time this approach degrades the control performance.
To reduce such errors, the variables having a large interaction should be in the same group,
and the variable interaction between different groups should be small.
Such a problem is called a \emph{graph partitioning problem}, which is known to be an NP-hard problem, but there are some approximation methods with adequate accuracy.
For the actual implementation, we used the \emph{Metis} software \reva{(Apache license 2.0)}, which is a
widely used solver for graph partitioning problems, to break up the large-scale problem into several groups having fewer than 64 variables~\cite{Karypis1998Fast}.
Figure~\ref{fig:metis} shows the result of the graph partitioning of the city of $ L = 50 $ into 42 groups using Metis,
where we certainly see that the adjacent intersections, i.e., the strongly interacting variables, are included in the same group.

\begin{figure}[t]
    \centering
    \includegraphics[width=80mm,bb=0.000000 0.000000 393.083371 267.635810]{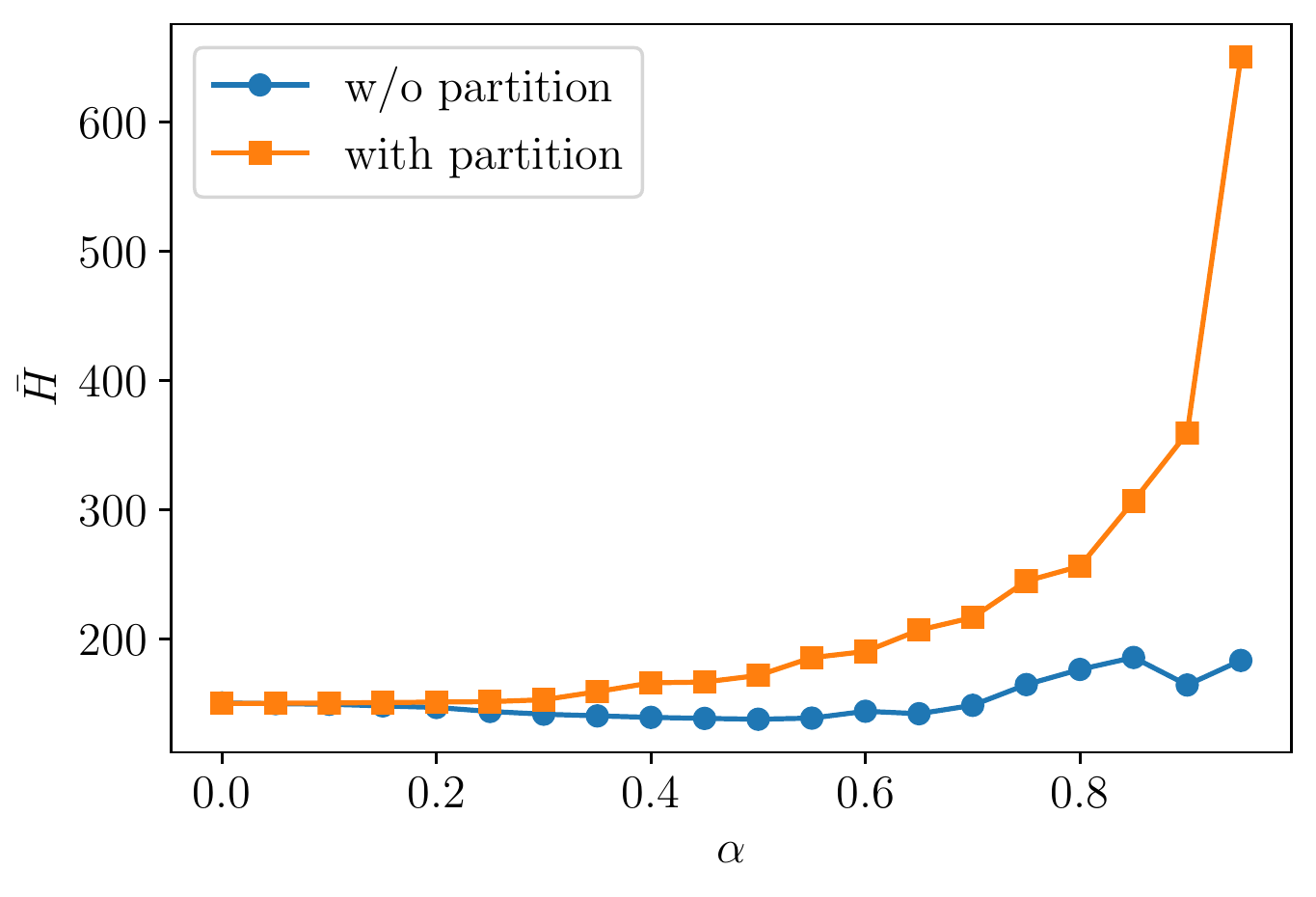}
    \caption{\rev{Time average of the objective function for the various $\alpha$, where $L=8$ and $\eta=1.0$. Orange-squares: with graph partitioning, blue-circles: without graph partitioning.}}
    \label{fig:error_partition}
\end{figure}

\rev{
 We evaluate the effect of this graph partitioning method.
 To this end, the simulated annealing optimization on a system with $L=8$ is performed, and results are compared between the no-partitioned and quadruple-partitioned cases.
  Figure~\ref{fig:error_partition} shows the time average of the objective function for various $\alpha$.
  The values with partitioning are larger than those without partitioning, where the difference between these values represents the error caused by partitioning.
  The error increases with larger $\alpha$, indicating that the larger the straight driving rate of vehicles, the more the partitioning has a negative effect.
In Fig.~\ref{fig:hamiltonian}, the quantum annealing has advantage at large $\alpha$, and thus a higher performance signal control should be achieved once the D-wave without partitioning is realized.
}


\begin{thebibliography}{10}
\urlstyle{rm}
\expandafter\ifx\csname url\endcsname\relax
  \def\url#1{\texttt{#1}}\fi
\expandafter\ifx\csname urlprefix\endcsname\relax\def\urlprefix{URL }\fi
\expandafter\ifx\csname doiprefix\endcsname\relax\def\doiprefix{DOI: }\fi
\providecommand{\bibinfo}[2]{#2}
\providecommand{\eprint}[2][]{\url{#2}}

\bibitem{Zhang2011DataDriven}
\bibinfo{author}{Zhang, J.} \emph{et~al.}
\newblock \bibinfo{journal}{\bibinfo{title}{Data-{{Driven Intelligent
  Transportation Systems}}: {{A Survey}}}}.
\newblock {\emph{\JournalTitle{IEEE Transactions on Intelligent Transportation
  Systems}}} \textbf{\bibinfo{volume}{12}}, \bibinfo{pages}{1624--1639},
  \doiprefix\url{10.1109/TITS.2011.2158001} (\bibinfo{year}{2011}).

\bibitem{Bishop2005Intelligent}
\bibinfo{author}{Bishop, R.}
\newblock \emph{\bibinfo{title}{Intelligent {{Vehicle Technology}} and
  {{Trends}}}} (\bibinfo{publisher}{{Artech House}}, \bibinfo{year}{2005}).

\bibitem{Cheng2015D2D}
\bibinfo{author}{Cheng, X.}, \bibinfo{author}{Yang, L.} \&
  \bibinfo{author}{Shen, X.}
\newblock \bibinfo{journal}{\bibinfo{title}{{{D2D}} for {{Intelligent
  Transportation Systems}}: {{A Feasibility Study}}}}.
\newblock {\emph{\JournalTitle{IEEE Transactions on Intelligent Transportation
  Systems}}} \textbf{\bibinfo{volume}{16}}, \bibinfo{pages}{1784--1793},
  \doiprefix\url{10.1109/TITS.2014.2377074} (\bibinfo{year}{2015}).

\bibitem{Papageorgiou2003Review}
\bibinfo{author}{Papageorgiou, M.}, \bibinfo{author}{Diakaki, C.},
  \bibinfo{author}{Dinopoulou, V.}, \bibinfo{author}{Kotsialos, A.} \&
  \bibinfo{author}{{Yibing Wang}}.
\newblock \bibinfo{journal}{\bibinfo{title}{Review of {{Road Traffic Control
  Strategies}}}}.
\newblock {\emph{\JournalTitle{Proceedings of the IEEE}}}
  \textbf{\bibinfo{volume}{91}}, \bibinfo{pages}{2043--2067},
  \doiprefix\url{10.1109/JPROC.2003.819610} (\bibinfo{year}{2003}).

\bibitem{Wei2019Survey}
\bibinfo{author}{Wei, H.}, \bibinfo{author}{Zheng, G.}, \bibinfo{author}{Gayah,
  V.} \& \bibinfo{author}{Li, Z.}
\newblock \bibinfo{journal}{\bibinfo{title}{A {{Survey}} on {{Traffic Signal
  Control Methods}}}}.
\newblock {\emph{\JournalTitle{arXiv:1904.08117 [cs, stat]}}}
  (\bibinfo{year}{2019}).
\newblock \eprint{1904.08117}.

\bibitem{Gokulan2010Distributed}
\rev{\bibinfo{author}{Gokulan, B.~P.} \& \bibinfo{author}{Srinivasan, D.}
\newblock \bibinfo{journal}{\bibinfo{title}{Distributed {{Geometric Fuzzy
  Multiagent Urban Traffic Signal Control}}}}.
\newblock {\emph{\JournalTitle{IEEE Transactions on Intelligent Transportation
  Systems}}} \textbf{\bibinfo{volume}{11}}, \bibinfo{pages}{714--727},
  \doiprefix\url{10.1109/TITS.2010.2050688} (\bibinfo{year}{2010}).}

\bibitem{Garcia-Nieto2012Swarm}
\rev{\bibinfo{author}{{Garc{\'i}a-Nieto}, J.}, \bibinfo{author}{Alba, E.} \&
  \bibinfo{author}{Carolina~Olivera, A.}
\newblock \bibinfo{journal}{\bibinfo{title}{Swarm intelligence for traffic
  light scheduling: {{Application}} to real urban areas}}.
\newblock {\emph{\JournalTitle{Engineering Applications of Artificial
  Intelligence}}} \textbf{\bibinfo{volume}{25}}, \bibinfo{pages}{274--283},
  \doiprefix\url{10.1016/j.engappai.2011.04.011} (\bibinfo{year}{2012}).}

\bibitem{Srinivasan2006Neural}
\rev{\bibinfo{author}{Srinivasan, D.}, \bibinfo{author}{Choy, M.} \&
  \bibinfo{author}{Cheu, R.}
\newblock \bibinfo{journal}{\bibinfo{title}{Neural {{Networks}} for
  {{Real}}-{{Time Traffic Signal Control}}}}.
\newblock {\emph{\JournalTitle{IEEE Transactions on Intelligent Transportation
  Systems}}} \textbf{\bibinfo{volume}{7}}, \bibinfo{pages}{261--272},
  \doiprefix\url{10.1109/TITS.2006.874716} (\bibinfo{year}{2006}).}

\bibitem{Arel2010Reinforcement}
\bibinfo{author}{Arel, I.}, \bibinfo{author}{Liu, C.},
  \bibinfo{author}{Urbanik, T.} \& \bibinfo{author}{Kohls, A.~G.}
\newblock \bibinfo{journal}{\bibinfo{title}{Reinforcement {{Learning}}-based
  {{Multi}}-agent {{System}} for {{Network Traffic Signal Control}}}}.
\newblock {\emph{\JournalTitle{IET Intelligent Transport Systems}}}
  \textbf{\bibinfo{volume}{4}}, \bibinfo{pages}{128--135},
  \doiprefix\url{10.1049/iet-its.2009.0070} (\bibinfo{year}{2010}).

\bibitem{Nishi2018Traffic}
\bibinfo{author}{Nishi, T.}, \bibinfo{author}{Otaki, K.},
  \bibinfo{author}{Hayakawa, K.} \& \bibinfo{author}{Yoshimura, T.}
\newblock \bibinfo{title}{Traffic {{Signal Control Based}} on {{Reinforcement
  Learning}} with {{Graph Convolutional Neural Nets}}}.
\newblock In \emph{\bibinfo{booktitle}{2018 21st International Conference on
  Intelligent Transportation Systems ({{ITSC}})}}, \bibinfo{pages}{877--883},
  \doiprefix\url{10.1109/ITSC.2018.8569301} (\bibinfo{organization}{{IEEE}},
  \bibinfo{year}{2018}).

\bibitem{Hunt1981Scoot}
\rev{\bibinfo{author}{Hunt, P.~B.}, \bibinfo{author}{Robertson, D.~I.},
  \bibinfo{author}{Bretherton, R.~D.} \& \bibinfo{author}{Winton, R.~I.}
\newblock \bibinfo{journal}{\bibinfo{title}{Scoot - {{A Traffic Responsive
  Method}} of {{Coordinating Signals}}}}.
\newblock {\emph{\JournalTitle{Publication of: Transport and Road Research
  Laboratory}}}  (\bibinfo{year}{1981}).}

\bibitem{Roess2004Traffic}
\bibinfo{author}{Roess, R.~P.}, \bibinfo{author}{Prassas, E.~S.} \&
  \bibinfo{author}{McShane, W.~R.}
\newblock \emph{\bibinfo{title}{Traffic Engineering}}
  (\bibinfo{publisher}{{Pearson/Prentice Hall}}, \bibinfo{year}{2004}).

\bibitem{Koonce2008Traffic}
\bibinfo{author}{Koonce, P.} \& \bibinfo{author}{Rodegerdts, L.}
\newblock \bibinfo{title}{Traffic {{Signal Timing Manual}}.}
\newblock \bibinfo{type}{Tech. Rep.}, \bibinfo{institution}{{United States.
  Federal Highway Administration}} (\bibinfo{year}{2008}).

\bibitem{Faouzi2011Data}
\rev{\bibinfo{author}{Faouzi, N.-E.~E.}, \bibinfo{author}{Leung, H.} \&
  \bibinfo{author}{Kurian, A.}
\newblock \bibinfo{journal}{\bibinfo{title}{Data {{Fusion}} in {{Intelligent
  Transportation Systems}}: {{Progress}} and {{Challenges}} \textendash{} {{A
  Survey}}}}.
\newblock {\emph{\JournalTitle{Information Fusion}}}
  \textbf{\bibinfo{volume}{12}}, \bibinfo{pages}{4--10},
  \doiprefix\url{10.1016/j.inffus.2010.06.001} (\bibinfo{year}{2011}).}

\bibitem{Khamis2012Multiobjective}
\bibinfo{author}{Khamis, M.~A.}, \bibinfo{author}{Gomaa, W.} \&
  \bibinfo{author}{{El-Shishiny}, H.}
\newblock \bibinfo{title}{Multi-objective {{Traffic Light Control System
  Based}} on {{Bayesian Probability Interpretation}}}.
\newblock In \emph{\bibinfo{booktitle}{2012 15th International {{IEEE}}
  Conference on Intelligent Transportation Systems}},
  \bibinfo{pages}{995--1000}, \doiprefix\url{10.1109/ITSC.2012.6338853}
  (\bibinfo{organization}{{IEEE}}, \bibinfo{year}{2012}).

\bibitem{Varaiya2013Maxpressure}
\bibinfo{author}{Varaiya, P.}
\newblock \bibinfo{title}{The {{Max}}-pressure {{Controller}} for {{Arbitrary
  Networks}} of {{Signalized Intersections}}}.
\newblock In \bibinfo{editor}{Ukkusuri, S.~V.} \& \bibinfo{editor}{Ozbay, K.}
  (eds.) \emph{\bibinfo{booktitle}{Advances in Dynamic Network Modeling in
  Complex Transportation Systems}}, \bibinfo{pages}{27--66},
  \doiprefix\url{10.1007/978-1-4614-6243-9\%0082}
  (\bibinfo{publisher}{{Springer New York}}, \bibinfo{address}{{New York, NY}},
  \bibinfo{year}{2013}).

\bibitem{Blum2003Metaheuristics}
\reva{
\bibinfo{author}{Blum, C.} \& \bibinfo{author}{Roli, A.}
\newblock \bibinfo{journal}{\bibinfo{title}{Metaheuristics in Combinatorial
  Optimization: {{Overview}} and Conceptual Comparison}}.
\newblock {\emph{\JournalTitle{ACM Computing Surveys}}}
  \textbf{\bibinfo{volume}{35}}, \bibinfo{pages}{268--308},
  \doiprefix\url{10.1145/937503.937505} (\bibinfo{year}{2003}).
}

\bibitem{Puchinger2005Combining}
\reva{
\bibinfo{author}{Puchinger, J.} \& \bibinfo{author}{Raidl, G.~R.}
\newblock \bibinfo{title}{Combining {{Metaheuristics}} and {{Exact Algorithms}}
  in {{Combinatorial Optimization}}: {{A Survey}} and {{Classification}}}.
\newblock In \bibinfo{editor}{Mira, J.} \& \bibinfo{editor}{{\'A}lvarez, J.~R.}
  (eds.) \emph{\bibinfo{booktitle}{Artificial {{Intelligence}} and {{Knowledge
  Engineering Applications}}: {{A Bioinspired Approach}}}}, Lecture {{Notes}}
  in {{Computer Science}}, \bibinfo{pages}{41--53},
  \doiprefix\url{10.1007/11499305_5} (\bibinfo{publisher}{{Springer}},
  \bibinfo{address}{{Berlin, Heidelberg}}, \bibinfo{year}{2005}).
}

\bibitem{Chakroun2013Combining}
\reva{
\bibinfo{author}{Chakroun, I.}, \bibinfo{author}{Melab, N.},
  \bibinfo{author}{Mezmaz, M.} \& \bibinfo{author}{Tuyttens, D.}
\newblock \bibinfo{journal}{\bibinfo{title}{Combining Multi-core and {{GPU}}
  Computing for Solving Combinatorial Optimization Problems}}.
\newblock {\emph{\JournalTitle{Journal of Parallel and Distributed Computing}}}
  \textbf{\bibinfo{volume}{73}}, \bibinfo{pages}{1563--1577},
  \doiprefix\url{10.1016/j.jpdc.2013.07.023} (\bibinfo{year}{2013}).
}

\bibitem{Inagaki2016coherent}
\rev{\bibinfo{author}{Inagaki, T.} \emph{et~al.}
\newblock \bibinfo{journal}{\bibinfo{title}{A Coherent {{Ising}} Machine for
  2000-node Optimization Problems}}.
\newblock {\emph{\JournalTitle{Science}}} \textbf{\bibinfo{volume}{354}},
  \bibinfo{pages}{603--606}, \doiprefix\url{10.1126/science.aah4243}
  (\bibinfo{year}{2016})}.

\bibitem{Hamerly2019Experimental}
\bibinfo{author}{Hamerly, R.} \emph{et~al.}
\newblock \bibinfo{journal}{\bibinfo{title}{Experimental {{Investigation}} of
  {{Performance Differences Between Coherent Ising Machines}} and a {{Quantum
  Annealer}}}}.
\newblock {\emph{\JournalTitle{Science Advances}}}
  \textbf{\bibinfo{volume}{5}}, \doiprefix\url{10.1126/sciadv.aau0823}
  (\bibinfo{year}{2019}).


\bibitem{Goto2019Combinatorial}
\reva{
\bibinfo{author}{Goto, H.}, \bibinfo{author}{Tatsumura, K.} \&
  \bibinfo{author}{Dixon, A.~R.}
\newblock \bibinfo{journal}{\bibinfo{title}{Combinatorial Optimization by
  Simulating Adiabatic Bifurcations in Nonlinear {{Hamiltonian}} Systems}}.
\newblock {\emph{\JournalTitle{Science Advances}}}
  \textbf{\bibinfo{volume}{5}}, \bibinfo{pages}{eaav2372},
  \doiprefix\url{10.1126/sciadv.aav2372} (\bibinfo{year}{2019}).
}

\bibitem{Matsubara2018IsingModel}
\bibinfo{author}{Matsubara, S.} \emph{et~al.}
\newblock \bibinfo{title}{Ising-{{Model Optimizer}} with {{Parallel}}-{{Trial
  Bit}}-{{Sieve Engine}}}.
\newblock In \bibinfo{editor}{Barolli, L.} \& \bibinfo{editor}{Terzo, O.}
  (eds.) \emph{\bibinfo{booktitle}{Complex, {{Intelligent}}, and {{Software
  Intensive Systems}}}}, Advances in {{Intelligent Systems}} and {{Computing}},
  \bibinfo{pages}{432--438}, \doiprefix\url{10.1007/978-3-319-61566-0_39}
  (\bibinfo{publisher}{{Springer International Publishing}},
  \bibinfo{address}{{Cham}}, \bibinfo{year}{2018}).


\bibitem{Aramon2019PhysicsInspired}
\reva{
\bibinfo{author}{Aramon, M.} \emph{et~al.}
\newblock \bibinfo{journal}{\bibinfo{title}{Physics-{{Inspired Optimization}}
  for {{Quadratic Unconstrained Problems Using}} a {{Digital Annealer}}}}.
\newblock {\emph{\JournalTitle{Frontiers in Physics}}}
  \textbf{\bibinfo{volume}{7}}, \doiprefix\url{10.3389/fphy.2019.00048}
  (\bibinfo{year}{2019}).
}

\bibitem{Kadowaki1998Quantum}
\bibinfo{author}{Kadowaki, T.} \& \bibinfo{author}{Nishimori, H.}
\newblock \bibinfo{journal}{\bibinfo{title}{Quantum {{Annealing}} in the
  {{Transverse Ising Model}}}}.
\newblock {\emph{\JournalTitle{Physical Review E}}}
  \textbf{\bibinfo{volume}{58}}, \bibinfo{pages}{5355--5363},
  \doiprefix\url{10.1103/PhysRevE.58.5355} (\bibinfo{year}{1998}).

\bibitem{Johnson2011Quantum}
\bibinfo{author}{Johnson, M.~W.} \emph{et~al.}
\newblock \bibinfo{journal}{\bibinfo{title}{Quantum {{Annealing}} with
  {{Manufactured Spins}}}}.
\newblock {\emph{\JournalTitle{Nature}}} \textbf{\bibinfo{volume}{473}},
  \bibinfo{pages}{194--198}, \doiprefix\url{10.1038/nature10012}
  (\bibinfo{year}{2011}).

\bibitem{Das2008Colloquium}
\reva{
\bibinfo{author}{Das, A.} \& \bibinfo{author}{Chakrabarti, B.~K.}
\newblock \bibinfo{journal}{\bibinfo{title}{Colloquium: {{Quantum}} Annealing
  and Analog Quantum Computation}}.
\newblock {\emph{\JournalTitle{Reviews of Modern Physics}}}
  \textbf{\bibinfo{volume}{80}}, \bibinfo{pages}{1061--1081},
  \doiprefix\url{10.1103/RevModPhys.80.1061} (\bibinfo{year}{2008}).
}

\bibitem{King2015Benchmarking}
\bibinfo{author}{King, J.}, \bibinfo{author}{Yarkoni, S.},
  \bibinfo{author}{Nevisi, M.~M.}, \bibinfo{author}{Hilton, J.~P.} \&
  \bibinfo{author}{McGeoch, C.~C.}
\newblock \bibinfo{journal}{\bibinfo{title}{Benchmarking a {{Quantum Annealing
  Processor}} with the {{Time}}-to-{{Target Metric}}}}.
\newblock {\emph{\JournalTitle{arXiv:1508.05087 [quant-ph]}}}
  (\bibinfo{year}{2015}).
\newblock \eprint{1508.05087}.

\bibitem{McGeoch2013Experimental}
\bibinfo{author}{McGeoch, C.~C.} \& \bibinfo{author}{Wang, C.}
\newblock \bibinfo{title}{Experimental {{Evaluation}} of an {{Adiabiatic
  Quantum System}} for {{Combinatorial Optimization}}}.
\newblock In \emph{\bibinfo{booktitle}{Proceedings of the {{ACM}} International
  Conference on Computing Frontiers}}, {{CF}} '13, \bibinfo{pages}{23},
  \doiprefix\url{10.1145/2482767.2482797}. \bibinfo{organization}{ACM}
  (\bibinfo{publisher}{{Association for Computing Machinery}},
  \bibinfo{address}{{New York, NY, USA}}, \bibinfo{year}{2013}).

\bibitem{Venturelli2016Quantum}
\bibinfo{author}{Venturelli, D.}, \bibinfo{author}{Marchand, D. J.~J.} \&
  \bibinfo{author}{Rojo, G.}
\newblock \bibinfo{journal}{\bibinfo{title}{Quantum {{Annealing
  Implementation}} of {{Job}}-{{Shop Scheduling}}}}.
\newblock {\emph{\JournalTitle{arXiv:1506.08479 [quant-ph]}}}
  (\bibinfo{year}{2016}).
\newblock \eprint{1506.08479}.

\bibitem{OMalley2018Nonnegative}
\bibinfo{author}{O'Malley, D.}, \bibinfo{author}{Vesselinov, V.~V.},
  \bibinfo{author}{Alexandrov, B.~S.} \& \bibinfo{author}{Alexandrov, L.~B.}
\newblock \bibinfo{journal}{\bibinfo{title}{Nonnegative/{{Binary Matrix
  Factorization}} with a {{D}}-{{Wave Quantum Annealer}}}}.
\newblock {\emph{\JournalTitle{PLOS ONE}}} \textbf{\bibinfo{volume}{13}},
  \bibinfo{pages}{e0206653}, \doiprefix\url{10.1371/journal.pone.0206653}
  (\bibinfo{year}{2018}).

\bibitem{Ohzeki2018Optimization}
\bibinfo{author}{Ohzeki, M.}, \bibinfo{author}{Okada, S.},
  \bibinfo{author}{Terabe, M.} \& \bibinfo{author}{Taguchi, S.}
\newblock \bibinfo{journal}{\bibinfo{title}{Optimization of {{Neural Networks
  Via Finite}}-{{Value Quantum Fluctuations}}}}.
\newblock {\emph{\JournalTitle{Scientific Reports}}}
  \textbf{\bibinfo{volume}{8}}, \bibinfo{pages}{1--10},
  \doiprefix\url{10.1038/s41598-018-28212-4} (\bibinfo{year}{2018}).

\bibitem{Inoue2020Model}
\bibinfo{author}{Inoue, D.} \& \bibinfo{author}{Yoshida, H.}
\newblock \bibinfo{journal}{\bibinfo{title}{Model {{Predictive Control}} for
  {{Finite Input Systems}} using the {{D}}-{{Wave Quantum Annealer}}}}.
\newblock {\emph{\JournalTitle{Scientific Reports}}}
  \textbf{\bibinfo{volume}{10}}, \bibinfo{pages}{1--10},
  \doiprefix\url{10.1038/s41598-020-58081-9} (\bibinfo{year}{2020}).

\bibitem{Ayanzadeh2020Reinforcementa}
\reva{
\bibinfo{author}{Ayanzadeh, R.}, \bibinfo{author}{Halem, M.} \&
  \bibinfo{author}{Finin, T.}
\newblock \bibinfo{journal}{\bibinfo{title}{Reinforcement {{Quantum
  Annealing}}: {{A Hybrid Quantum Learning Automata}}}}.
\newblock {\emph{\JournalTitle{Scientific Reports}}}
  \textbf{\bibinfo{volume}{10}}, \bibinfo{pages}{7952},
  \doiprefix\url{10.1038/s41598-020-64078-1} (\bibinfo{year}{2020}).
}

\bibitem{Yang1952Spontaneous}
\bibinfo{author}{Yang, C.~N.}
\newblock \bibinfo{journal}{\bibinfo{title}{The {{Spontaneous Magnetization}}
  of a {{Two}}-{{Dimensional Ising Model}}}}.
\newblock {\emph{\JournalTitle{Physical Review}}}
  \textbf{\bibinfo{volume}{85}}, \bibinfo{pages}{808--816},
  \doiprefix\url{10.1103/PhysRev.85.808} (\bibinfo{year}{1952}).

\bibitem{McCoy2014twodimensional}
\bibinfo{author}{McCoy, B.~M.} \& \bibinfo{author}{Wu, T.~T.}
\newblock \emph{\bibinfo{title}{The Two-Dimensional {{Ising}} Model}}
  (\bibinfo{publisher}{{Courier Corporation}}, \bibinfo{year}{2014}).

\bibitem{Binder1981Finite}
\bibinfo{author}{Binder, K.}
\newblock \bibinfo{journal}{\bibinfo{title}{Finite {{Size Scaling Analysis}} of
  {{Ising Model Block Distribution Functions}}}}.
\newblock {\emph{\JournalTitle{Zeitschrift f\"ur Physik B Condensed Matter}}}
  \textbf{\bibinfo{volume}{43}}, \bibinfo{pages}{119--140},
  \doiprefix\url{10.1007/BF01293604} (\bibinfo{year}{1981}).

\bibitem{Glauber1963Time}
\bibinfo{author}{Glauber, R.~J.}
\newblock \bibinfo{journal}{\bibinfo{title}{Time-{{Dependent Statistics}} of
  the {{Ising Model}}}}.
\newblock {\emph{\JournalTitle{Journal of Mathematical Physics}}}
  \textbf{\bibinfo{volume}{4}}, \bibinfo{pages}{294--307},
  \doiprefix\url{10.1063/1.1703954} (\bibinfo{year}{1963}).

\bibitem{Suzuki2013Chaotic}
\bibinfo{author}{Suzuki, H.}, \bibinfo{author}{Imura, J.-i.} \&
  \bibinfo{author}{Aihara, K.}
\newblock \bibinfo{journal}{\bibinfo{title}{Chaotic {{Ising}}-{{Like Dynamics}}
  in {{Traffic Signals}}}}.
\newblock {\emph{\JournalTitle{Scientific Reports}}}
  \textbf{\bibinfo{volume}{3}}, \bibinfo{pages}{1--6},
  \doiprefix\url{10.1038/srep01127} (\bibinfo{year}{2013}).

\bibitem{Suman2006survey}
\bibinfo{author}{Suman, B.} \& \bibinfo{author}{Kumar, P.}
\newblock \bibinfo{journal}{\bibinfo{title}{A {{Survey}} of {{Simulated
  Annealing}} as a {{Tool}} for {{Single}} and {{Multiobjective
  Optimization}}}}.
\newblock {\emph{\JournalTitle{Journal of the Operational Research Society}}}
  \textbf{\bibinfo{volume}{57}}, \bibinfo{pages}{1143--1160},
  \doiprefix\url{10.1057/palgrave.jors.2602068} (\bibinfo{year}{2006}).

\bibitem{Boothby2016Fast}
\bibinfo{author}{Boothby, T.}, \bibinfo{author}{King, A.~D.} \&
  \bibinfo{author}{Roy, A.}
\newblock \bibinfo{journal}{\bibinfo{title}{Fast {{Clique Minor Generation}} in
  {{Chimera Qubit Connectivity Graphs}}}}.
\newblock {\emph{\JournalTitle{Quantum Information Processing}}}
  \textbf{\bibinfo{volume}{15}}, \bibinfo{pages}{495--508},
  \doiprefix\url{10.1007/s11128-015-1150-6} (\bibinfo{year}{2016}).

\bibitem{W.Johnson2018Future}
\rev{\bibinfo{author}{Johnson, M.~W.}
\newblock \bibinfo{title}{Future {{Hardware Directions}} of {{Quantum
  Annealing}}}.
\newblock In \emph{\bibinfo{booktitle}{Qubits {{Europe}} 2018 {{D}}-{{Wave
  Users Conference}}}} (\bibinfo{address}{{Munich}}, \bibinfo{year}{2018}).}

\bibitem{Ayanzadeh2020PostQuantum}
\reva{
\bibinfo{author}{Ayanzadeh, R.}, \bibinfo{author}{Dorband, J.},
  \bibinfo{author}{Halem, M.} \& \bibinfo{author}{Finin, T.}
\newblock \bibinfo{journal}{\bibinfo{title}{Post-{{Quantum
  Error}}-{{Correction}} for {{Quantum Annealers}}}}.
\newblock {\emph{\JournalTitle{arXiv:2010.00115 [quant-ph]}}}
  (\bibinfo{year}{2020}).
\newblock \eprint{2010.00115}.
}

\bibitem{manual}
\reva{
 \bibinfo {howpublished} {See \url{https://docs.dwavesys.com/docs/latest/c\_solver\_2.htm} for the VFYC solver.}
}

\bibitem{Cai2014practicala}
\reva{
\bibinfo{author}{Cai, J.}, \bibinfo{author}{Macready, W.~G.} \&
  \bibinfo{author}{Roy, A.}
\newblock \bibinfo{journal}{\bibinfo{title}{A practical heuristic for finding
  graph minors}}.
\newblock {\emph{\JournalTitle{arXiv:1406.2741 [quant-ph]}}}
  (\bibinfo{year}{2014}).
\newblock \eprint{1406.2741}.
}

\bibitem{Karypis1998Fast}
\bibinfo{author}{Karypis, G.} \& \bibinfo{author}{Kumar, V.}
\newblock \bibinfo{journal}{\bibinfo{title}{A {{Fast}} and {{High Quality
  Multilevel Scheme}} for {{Partitioning Irregular Graphs}}}}.
\newblock {\emph{\JournalTitle{SIAM Journal on Scientific Computing}}}
  \textbf{\bibinfo{volume}{20}}, \bibinfo{pages}{359--392}
  (\bibinfo{year}{1998}).






\end{thebibliography}

\section*{Acknowledgements}
The authors would like to thank Dr.~Kiyosumi Kidono of Toyota Central R\&D Labs. for the useful discussions.
This research was performed as a project of Intelligent Mobility Society Design, Social Cooperation Program (Next Generation Artificial Intelligence Research Center, the University of Tokyo with Toyota Central R\&D Labs., Inc.).

\section*{Author contributions statement}
D.I. conceived and developed the concept and carried out all the experiments. 
D.I. and H.Y. analyzed the results and wrote the manuscript.
A.O., T.M., K.A., and H.Y. designed the research plan and reviewed the manuscript.

\section*{Additional information}

\textbf{Competing interests}: The authors declare no competing financial interests.

\end{document}